\documentclass[11pt]{amsart}
\usepackage{amssymb}
\usepackage{graphicx}

\newtheorem{lemma}{Lemma}[section]
\newtheorem{theorem}{Theorem}[section]
\newtheorem{proposition}{Proposition}[section]
\newcommand{\adfQED}{\hfill $\square$}
\newcommand{\adfmod}[1]{~(\mathrm{mod}~#1)}

\begin{document}
\title{Designs for graphs with six vertices and ten edges - II}
\author{A. D. Forbes}
\address{School of Mathematics and Statistics\\
        The Open University\\
        Walton Hall, Milton Keynes MK7 6AA, UK}
\email{anthony.d.forbes@gmail.com}
\author{T. S. Griggs}
\address{School of Mathematics and Statistics\\
        The Open University\\
        Walton Hall, Milton Keynes MK7 6AA, UK}
\email{terry.griggs@open.ac.uk}
\date{19 Apr 2020} 

\maketitle
\begin{abstract}
The design spectrum has been determined for ten of the 15 graphs with six vertices and ten edges.
In this paper we solve the design spectrum problem for the remaining five graphs with three possible exceptions.
\end{abstract}

\section{Introduction}\label{sec:Introduction}
Throughout this paper all graphs are simple.
Let $G$ be a graph. If the edge set of a graph $K$ can be partitioned into edge sets
of graphs each isomorphic to $G$, we say that there exists a {\em decomposition} of $K$ into $G$.
In the case where $K$ is the complete graph $K_n$ we refer to the decomposition as a
$G$-{\em design} of order $n$. The {\em design spectrum} of $G$ is the set of non-negative integers $n$ for
which there exists a $G$-design of order $n$.
We refer the reader to the survey article of Adams, Bryant and Buchanan, \cite{AdamsBryantBuchanan2008} and,
for more up to date results, the Web site maintained by Bryant and McCourt, \cite{BryantMcCourt9999}.
If the graph $G$ has $v$ vertices, $e$ edges, and if $d$ is the greatest common divisor of the vertex degrees,
then a $G$-design of order $n$ can exist only if the following conditions hold:
\begin{equation}\label{eqn:admissible}
\left\{\begin{array}{rl}
\textrm{(i)}   & n \le 1 \textrm{~or~} n \ge v,\\
\textrm{(ii)}  &  n-1 \equiv 0 \adfmod{d},\\
\textrm{(iii)} & n(n-1) \equiv 0 \adfmod{2 e}.
\end{array}\right.
\end{equation}
Except where (i) of (\ref{eqn:admissible}) applies, adding an isolated vertex to a graph does not affect
its design spectrum.

The problem for small graphs has attracted attention.
As far as the authors are aware, the design spectrum problem has been completely solved for
all graphs with up to five vertices and all graphs with six vertices and up to 9 edges.
For details and references, see \cite{AdamsBryantBuchanan2008} and \cite{BryantMcCourt9999}, and for more recent results,
\cite{GeLing2008}, 
\cite{Kolotoglu2013}, 
\cite{GeHuKolotogluWei2015}, 
\cite{ForbesGriggsForbes2017}, 
\cite{AdamsBillingtonHoffman1997}, 
\cite{GuyBeineke1968}, 
\cite{MullinPoploveZhu1987}, 
\cite{KangZhaoMa2008}, 
\cite{ForbesGriggs2017}, 
\cite{Kolotoglu2019}. 
In Table~\ref{tab:graphs} we list the 15 graphs with six vertices and ten edges.
The numbering in the first column corresponds to the ordering of the ten-edge graphs within the
list of all 156 graphs with six vertices, available at \cite{McKay9999}.
The second column identifies the graphs as they appear in {\em An Atlas of Graphs} by Read \& Wilson, \cite{ReadWilson1998}.
The third column contains the edge sets which we use in our computations; the vertices are labelled in non-increasing order of degree.
The spectrum has been completely determined for graphs
$n_9$ ($K_5$ with an additional, isolated vertex), \cite{Hanani1961},
$n_{11}$, \cite{AdamsBryantKhodkar2000}, as well as
$n_1$, $n_2$, $n_4$, $n_5$, $n_{7}$, $n_{12}$, $n_{14}$ and $n_{15}$, \cite{ForbesGriggsForbes2019}.
The purpose of this paper is to address the remaining five graphs with six vertices and ten edges:
$$n_{3},~ n_{6},~ n_{8},~ n_{10} \textrm{~~and~~} n_{13}.$$
\begin{table}[t!] 
\caption{The 15 graphs with 6 vertices and 10 edges}
\label{tab:graphs}
\begin{center}
\begin{tabular}{@{}l@{~~}l@{~~}l@{}}
  $n_{1 }$ & G179 & \{\{4,3\},\{4,2\},\{4,1\},\{6,2\},\{6,1\},\{5,2\},\{5,1\},\{3,2\},\{3,1\},\{2,1\}\}\\ 
  $n_{2 }$ & G180 & \{\{4,3\},\{4,2\},\{4,1\},\{6,3\},\{6,1\},\{5,2\},\{5,1\},\{3,2\},\{3,1\},\{2,1\}\}\\ 
  $n_{3 }$ & G177 & \{\{5,3\},\{5,2\},\{5,1\},\{4,3\},\{4,2\},\{4,1\},\{6,1\},\{3,2\},\{3,1\},\{2,1\}\}\\ 
  $n_{4 }$ & G182 & \{\{5,3\},\{5,2\},\{5,1\},\{4,3\},\{4,2\},\{4,1\},\{6,2\},\{6,1\},\{3,1\},\{2,1\}\}\\ 
  $n_{5 }$ & G186 & \{\{5,3\},\{5,2\},\{5,1\},\{4,3\},\{4,2\},\{4,1\},\{6,3\},\{6,2\},\{3,1\},\{2,1\}\}\\ 
  $n_{6 }$ & G189 & \{\{6,2\},\{6,3\},\{6,1\},\{5,2\},\{5,3\},\{5,1\},\{4,2\},\{4,3\},\{4,1\},\{2,1\}\}\\ 
  $n_{7 }$ & G183 & \{\{5,3\},\{5,2\},\{5,1\},\{4,6\},\{4,2\},\{4,1\},\{3,2\},\{3,1\},\{6,1\},\{2,1\}\}\\ 
  $n_{8 }$ & G190 & \{\{6,4\},\{6,2\},\{6,1\},\{5,3\},\{5,2\},\{5,1\},\{4,2\},\{4,1\},\{3,2\},\{3,1\}\}\\ 
  $n_{9 }$ & G176 & \{\{5,4\},\{5,3\},\{5,2\},\{5,1\},\{4,3\},\{4,2\},\{4,1\},\{3,2\},\{3,1\},\{2,1\}\}\\ 
  $n_{10}$ & G178 & \{\{4,3\},\{4,2\},\{4,5\},\{4,1\},\{6,1\},\{3,2\},\{3,5\},\{3,1\},\{2,5\},\{2,1\}\}\\ 
  $n_{11}$ & G181 & \{\{4,3\},\{4,5\},\{4,2\},\{4,1\},\{6,2\},\{6,1\},\{3,5\},\{3,2\},\{3,1\},\{2,1\}\}\\ 
  $n_{12}$ & G185 & \{\{3,2\},\{3,5\},\{3,4\},\{3,1\},\{6,4\},\{6,1\},\{2,5\},\{2,4\},\{2,1\},\{5,1\}\}\\ 
  $n_{13}$ & G187 & \{\{6,4\},\{6,3\},\{6,1\},\{5,3\},\{5,2\},\{5,1\},\{4,2\},\{4,1\},\{3,1\},\{2,1\}\}\\ 
  $n_{14}$ & G184 & \{\{3,6\},\{3,4\},\{3,2\},\{3,1\},\{5,4\},\{5,2\},\{5,1\},\{6,2\},\{4,1\},\{2,1\}\}\\ 
  $n_{15}$ & G188 & \{\{2,5\},\{2,4\},\{2,3\},\{2,1\},\{6,4\},\{6,3\},\{6,1\},\{5,3\},\{5,1\},\{4,1\}\}\\ 
\end{tabular}
\end{center}
\end{table}
\begin{figure}[t!]
\caption{Graphs with 6 vertices and 10 edges}
\label{fig:graph-pictures}
\begin{center}
\includegraphics[width=120mm]{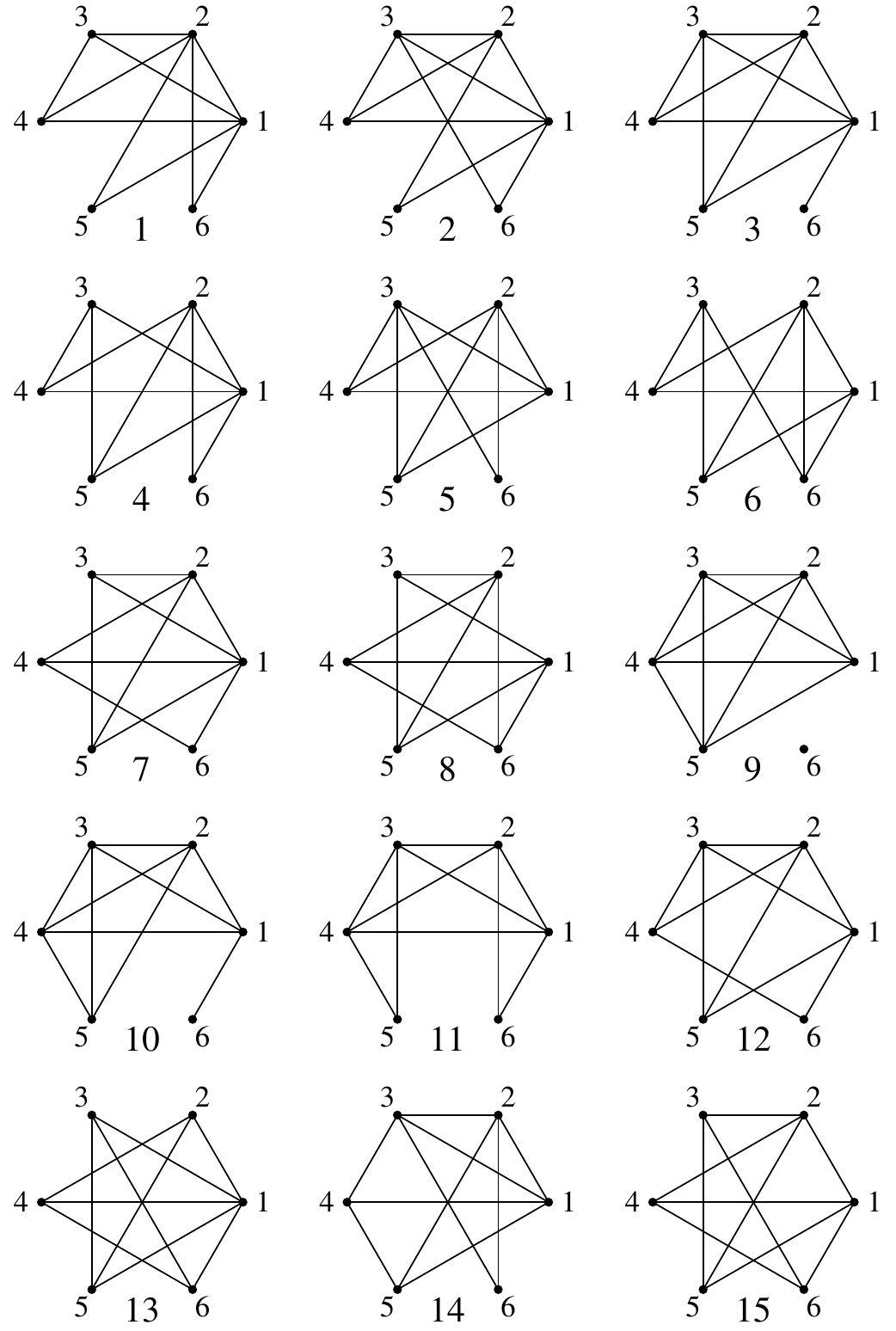}
\end{center}
\end{figure}
We prove the following.

\begin{theorem}
\label{thm:v6e10main}
$\mathrm{(i)}$ Designs of order $n$ exist for graphs $n_{3}$, $n_{6}$ and $n_{10}$ if and only if
$n \equiv 0,~ 1,~ 5,~ 16$ $\adfmod{20}$, except for $n = 5$ and except possibly for $n = 16$.

$\mathrm{(ii)}$  Designs of order $n$ exist for graph $n_{8}$ if and only if
$n \equiv 0,~ 1,~ 5,~ 16$ $\adfmod{20}$, except for $n \in \{5, 16\}$.

$\mathrm{(iii)}$ Designs of order $n$ exist for graph $n_{13}$ if and only if
$n \equiv 0,~ 1,~ 5,~ 16$ $\adfmod{20}$, except for $n \in \{5, 16, 20\}$.
\end{theorem}

Theorem~\ref{thm:v6e10main} is proved in Section~\ref{sec:v6e10main}.
For our computations and in the presentation of our results, we represent the labelled graph $n_i$ by a
subscripted ordered $6$-tuple $(z_1, z_2, \dots, z_6)_{i}$,
where $z_1 = 1$, $z_2 = 2$, \dots, $z_6 = 6$ give the edge sets in Table~\ref{tab:graphs} and the
illustrations in Figure~\ref{fig:graph-pictures}.

\section{Nonexistence results}\label{sec:v6e10nonexistence}

\begin{lemma}
\label{lem:nonexistence order 16}
A design of order $16$ does not exist for graph $n_{8}$ or graph $n_{13}$.
\end{lemma}
\begin{proof}
In the design there are 12 graphs whose vertices are labelled with symbols representing the vertices of $K_{16}$.
A vertex of $K_{16}$ has degree 15.

For $n_{8}$, the vertices of $n_{8}$ have degrees $(4,4,3,3,3,3)$.
A label must belong to one of two sets:
$A$: labels attached to vertices of degrees \{4,4,4,3\}, and
$B$: labels attached to vertices of degrees \{3,3,3,3,3\},
in each case the degrees sum to 15, the degree of a vertex of $K_{16}$.
Then $3|A| = 24$, $|A| + 5|B| = 48$ and hence $|A| = |B| = 8$.
Now consider pairs of labels of vertices of $K_{16}$ which are both in set $B$.
There are ${8 \choose 2} = 28$ of these.
But in $n_8$ only two pairs of vertices of degree 3, $\{3,5\}$ and $\{4,6\}$, are adjacent,
which creates only $12 \cdot 2 = 24$ pairs in total.

The vertices of $n_{13}$ have degrees $(5,3,3,3,3,3)$.
The only partitions of 15 into elements from $\{3,5\}$ are $\{5,5,5\}$ and $\{3,3,3,3,3\}$.
Let $A$ denote the set of labels that are attached to vertices of degrees $\{5,5,5\}$.
Then $3|A| = 12$ and hence $|A| = 4 > 1$. However, the design cannot have any $A$-$A$ pairs.
\end{proof}

\begin{lemma}
\label{lem:nonexistence order 20}
A design of order $20$ does not exist for graph $n_{13}$.
\end{lemma}
\begin{proof}
In the design there are 19 graphs whose vertices are labelled with symbols representing the vertices of $K_{20}$.
A vertex of $K_{20}$ has degree 19.

The vertices of $n_{13}$ have degrees $(5,3,3,3,3,3)$.
Since the only partition of 19 into elements from $\{3,5\}$ is $\{5,5,3,3,3\}$,
each label must be attached to two $n_{13}$ vertices of degree 5.
This is impossible.
\end{proof}

\section{Proof of Theorem~\ref{thm:v6e10main}}\label{sec:v6e10main}

We use Wilson's construction involving group divisible designs.
For this paper, a $K$-GDD of type $g_1^{t_1} \dots\ g_r^{t_r}$ is an ordered triple ($V,\mathcal{G},\mathcal{B}$)
where $V$ is a base set of cardinality $v = t_1 g_1 + \dots + t_r g_r$, $\mathcal{G}$ is a partition of
$V$ into $t_i$ subsets of cardinality $g_i$, $i = 1, \dots, r$, called \textit{groups} and $\mathcal{B}$
is a collection of subsets of cardinalities $k \in K$, called \textit{blocks}, which collectively have the
property that each pair of elements from different groups occurs in precisely one block but no pair of
elements from the same group occurs at all. A $\{k\}$-GDD is also called a $k$-GDD.
As is well known, if there exist $k-2$ mutually orthogonal Latin squares of side $q$, then there exists a $k$-GDD of type $q^k$.
So when $q$ is a prime power there exists a $q$-GDD of type $q^q$ and a $(q+1)$-GDD of type $q^{q+1}$
(obtained from affine and projective planes of order $q$ respectively).
A {\em parallel class} in a group divisible design is a subset of
the block set in which each element of the base set appears exactly once.
A $k$-GDD is called {\em resolvable}, and denoted by $k$-RGDD, if the entire set of blocks
can be partitioned into parallel classes.

Our primary construction is exactly the same as the one used in \cite{ForbesGriggsForbes2019}. We repeat it here for convenience.

\begin{proposition}
\label{prop:construction-main}
Let $i$, $t$, $p$, $q$ be positive integers.
Let $w$, $x$, $y$ be non-negative integers such that $x + y = w$ and $w \le 4t$.
Let $e = 0$ or $1$.
Suppose there exist decompositions into the graph $G$ of the complete graphs $K_{4i + e}$ and $K_{xp + yq + e}$ as well
as the complete multipartite graphs $K_{i,i,i,i}$, $K_{i,i,i,i,p}$ and $K_{i,i,i,i,q}$.
Then there exists a $G$-design of order $12it + 4i + xp + yq + e$.
\end{proposition}
\begin{proof}
See \cite[Proposition 2.1]{ForbesGriggsForbes2019}.
\end{proof}


\begin{sloppypar}
Before applying Proposition~\ref{prop:construction-main} we establish the existence of the various decompositions that we need
to make the construction work. With $i = 10$ and $p$, $q$ chosen from $\{10, 15, 20\}$, we require the following:
\begin{enumerate}
\item[]decompositions of $K_{10,10,10,10}$, $K_{10,10,10,10,10}$, $K_{10,10,10,10,15}$ and $K_{10,10,10,10,20}$;
\item[]design orders 40, 60,  80, 100, 140,  21, 41,  61,  81, 101,  25, 45,  65,  85, 105, 125,  36, 56,  76,  96, 116 and 136.
\end{enumerate}
Observe that the required design orders $xp + yq + e$ of Proposition~\ref{prop:construction-main} correspond to
the third term of the sums in the column headed `order' in Table~\ref{tab:main construction}.
The existence of these designs is proved in Lemmas~\ref{lem:basic designs direct} and \ref{lem:basic designs constructed}, but first,
in Lemma~\ref{lem:basic multipartite direct}, we give the decompositions of complete multipartite graphs that we will need for all of our constructions.
\end{sloppypar}

\begin{lemma}
\label{lem:basic multipartite direct}
\begin{sloppypar}
There exist decompositions of
$K_{10,10,10,10}$, $K_{15,15,15,15}$,
$K_{20,20,20,20}$, $K_{25,25,25,25}$,
$K_{10,10,10,15}$,
$K_{5,5,5,5,5}$, $K_{6,6,6,6,6}$, $K_{10,10,10,10,10}$,
$K_{10,10,10,10,15}$, $K_{10,10,10,10,20}$ and
$K_{4^6}$
into each of the five graphs.
\end{sloppypar}

There exist decompositions of $K_{10,10,10}$ and $K_{5,5,5,9}$ into graphs $n_{6}$ and $n_{8}$.

There exist decompositions of $K_{3,3,3,3,3}$ into graphs $n_{6}$ and $n_{10}$.

There exists a decomposition of $K_{21,21,21,21,21}$ into graph $n_{3}$.

There exist decompositions of $K_{8,8,8,8,8}$, $K_{8,8,8,8,3}$,
$K_{21,21,21,21,36}$, $K_{4^6,10}$, $K_{4^6,15}$,
$K_{1^{39},21}$, $K_{1^{55},25}$ and $K_{1^{99},41}$ into graph $n_{13}$.

There exist decompositions of $K_{4,4,4,4,4,7}$ and $K_{4^6,5}$ into graphs $n_{3}$, $n_{6}$, $n_{8}$ and $n_{10}$.
\end{lemma}
\begin{proof}
The decompositions are presented in Appendix~\ref{sec:App-M}.
\end{proof}

\begin{lemma}
\label{lem:basic designs direct}
Designs of orders $40$, $21$, $41$, $25$, $45$, $65$, $85$, $36$, $56$, $76$ and $116$ exist for all five graphs.

Designs of order $20$ exist for graphs $n_{3}$, $n_{6}$, $n_{8}$ and $n_{10}$.

Designs of orders $60$ and $80$ exist for graphs $n_{3}$, and $n_{10}$.

Designs of order $61$ exist for graphs $n_{3}$, $n_{10}$ and $n_{13}$.

Designs of order $105$ exist for graphs $n_{8}$ and $n_{13}$.

Designs of order $96$ exist for graphs $n_{3}$, $n_{10}$ and $n_{13}$.

A design of order $136$ exists for graph $n_{13}$.

Designs of order $156$ exist for graphs $n_{3}$, $n_{10}$ and $n_{13}$.
\end{lemma}
\begin{proof}
The decompositions are presented in Appendix~\ref{sec:App-D}.
\end{proof}

\begin{lemma}
\label{lem:basic designs constructed}
Designs of orders $100$, $140$,
$81$, $101$ and
$125$ exist for all five graphs.

Designs of orders $60$ and $80$ exist for graphs $n_{6}$, $n_{8}$ and $n_{13}$.

Designs of order $61$ exist for graphs $n_{6}$ and $n_{8}$.

Designs of order $105$ exist for graphs $n_{3}$, $n_{6}$ and $n_{10}$.

Designs of order $96$ exist for graphs $n_{6}$ and $n_{8}$.

Designs of order $136$ exist for graphs $n_{3}$, $n_{6}$, $n_{8}$ and $n_{10}$.
\end{lemma}
\begin{proof}
These designs are constructed.
We give only brief details by specifying the ingredients for Wilson's construction, namely
the complete graphs, the complete multipartite graphs and the group divisible designs.
It should be clear how the points of the GDDs are inflated and which GDDs are augmented by an extra point.
Decompositions of the ingredients exist by
Lemmas~\ref{lem:basic multipartite direct} and \ref{lem:basic designs direct}.


\begin{sloppypar}
For order 60 for graphs $n_{6}$ and $n_{8}$, use decompositions of
$K_{20}$ and $K_{10,10,10}$ with a 3-GDD of type $2^3$ (obtained from a Latin square of side 2).

For order 60 for graph $n_{13}$, use decompositions of
$K_{21}$ and $K_{1^{39},21}$.

For order 80 for graphs $n_{6}$ and $n_{8}$, use decompositions of
$K_{20}$, and $K_{20,20,20,20}$.

For order 80 for graph $n_{13}$, use decompositions of
$K_{25}$ and $K_{1^{55},25}$.

For order 100 for all five graphs, use decompositions of
$K_{25}$ and $K_{25,25,25,25}$.

For order 140 for graphs $n_{3}$, $n_{6}$, $n_{8}$ and $n_{10}$, use decompositions of
$K_{20}$ and $K_{10,10,10,10}$ with a 4-GDD of type $2^7$, \cite{BrouwerSchrijverHanani1977}, \cite{Ge2007}.

For order 140 for graph $n_{13}$, use decompositions of
$K_{41}$ and $K_{1^{99},41}$.


For order 61 for $n_{6}$ and $n_{8}$, use decompositions of
$K_{21}$ and $K_{10,10,10}$ with a 3-GDD of type $2^3$.

For order 81 for all five graphs, use decompositions of
$K_{21}$, and $K_{20,20,20,20}$.

For order 101 for all five graphs, use decompositions of
$K_{21}$ and $K_{5,5,5,5,5}$ with a 5-GDD of type $4^5$ (obtained from a projective plane of order 4).


For order 105 for graph $n_{3}$, use decompositions of
$K_{21}$ and $K_{21,21,21,21,21}$.

For order 105 for graphs $n_{6}$ and $n_{10}$, use decompositions of
$K_{21}$ and $K_{3,3,3,3,3}$ with a 5-GDD of type $7^5$ (obtained from 3 MOLS of side 7).

For order 125 for all five graphs, use decompositions of
$K_{25}$ and $K_{5,5,5,5,5}$ with a 5-GDD of type $5^5$.


For order 96 for graphs $n_{6}$ and $n_{8}$, use decompositions of
$K_{20}$, $K_{36}$ and $K_{5,5,5,9}$ with a 4-GDD of type $4^4$ (obtained from an affine plane of order 4).
\end{sloppypar}

For order 136 for graphs $n_{3}$, $n_{6}$, $n_{8}$ and $n_{10}$, use decompositions of
$K_{21}$, $K_{36}$ and $K_{4,4,4,4,4,7}$ with a 6-GDD of type $5^6$ (obtained from a projective plane of order 5).
\end{proof}

We are now ready to prove Theorem~\ref{thm:v6e10main}. For the main construction, we use Proposition~\ref{prop:construction-main} with $i = 10$,
and $p, q \in \{10, 15, 20\}$ as indicated in Table~\ref{tab:main construction}.
The inflated blocks of the group divisible design become multipartite graphs
$K_{10,10,10,10}$, $K_{10,10,10,10,10}$, $K_{10,10,10,10,15}$ and $K_{10,10,10,10,20}$.
With the decompositions of Lemmas~\ref{lem:basic multipartite direct}, \ref{lem:basic designs direct} and \ref{lem:basic designs constructed}
we obtain the designs listed in Table~\ref{tab:main construction}. Combining the residue classes modulo 120, we see that
we can construct designs of order $n$, $n \equiv 0$, 1, 5 and 16 (modulo 20) except for those orders listed under `missing values'.
\begin{table}[t!]
\caption{The main construction}
\label{tab:main construction}
\begin{center}
\begin{tabular}{@{}lcccccccl@{}}
 order             & $t$       & $w$ & $x$ & $p$ & $y$ & $q$ & $e$ & missing values\\
\hline
$120t + 40      $ & $t \ge 1$ & 0 & 0 &  - & 0 &  - & 0 & 40 \\
$120t + 40 + 140$ & $t \ge 2$ & 7 & 0 &  - & 7 & 20 & 0 & 60,\, 180,\, 300 \\
$120t + 40 +  40$ & $t \ge 1$ & 2 & 0 &  - & 2 & 20 & 0 & 80 \\
$120t + 40 +  60$ & $t \ge 1$ & 3 & 0 &  - & 3 & 20 & 0 & 100 \\
$120t + 40 +  80$ & $t \ge 1$ & 4 & 0 &  - & 4 & 20 & 0 & 120 \\
$120t + 40 + 100$ & $t \ge 2$ & 5 & 0 &  - & 5 & 20 & 0 & 20,\, 140,\, 260\\[1mm]
$120t + 40 +   1$ & $t \ge 1$ & 0 & 0 &  - & 0 &  - & 1 & 41 \\
$120t + 40 +  21$ & $t \ge 1$ & 1 & 0 &  - & 1 & 20 & 1 & 61 \\
$120t + 40 +  41$ & $t \ge 1$ & 2 & 0 &  - & 2 & 20 & 1 & 81 \\
$120t + 40 +  61$ & $t \ge 1$ & 3 & 0 &  - & 3 & 20 & 1 & 101 \\
$120t + 40 +  81$ & $t \ge 1$ & 4 & 0 &  - & 4 & 20 & 1 & 121 \\
$120t + 40 + 101$ & $t \ge 2$ & 5 & 0 &  - & 5 & 20 & 1 & 21,\, 141,\, 261\\[1mm]
$120t + 40 + 125$ & $t \ge 2$ & 7 & 3 & 15 & 4 & 20 & 0 & 45,\, 165,\, 285 \\
$120t + 40 +  25$ & $t \ge 1$ & 2 & 1 & 10 & 1 & 15 & 0 & 65 \\
$120t + 40 +  45$ & $t \ge 1$ & 3 & 3 & 15 & 0 &  - & 0 & 85 \\
$120t + 40 +  65$ & $t \ge 1$ & 4 & 3 & 15 & 1 & 20 & 0 & 105 \\
$120t + 40 +  85$ & $t \ge 2$ & 5 & 3 & 15 & 2 & 20 & 0 & 5,\, 125,\, 245 \\
$120t + 40 + 105$ & $t \ge 2$ & 6 & 3 & 15 & 3 & 20 & 0 & 25,\, 145,\, 265\\[1mm]
$120t + 40 + 136$ & $t \ge 2$ & 7 & 1 & 15 & 6 & 20 & 1 & 56,\, 176,\, 296 \\
$120t + 40 +  36$ & $t \ge 1$ & 2 & 1 & 15 & 1 & 20 & 1 & 76 \\
$120t + 40 +  56$ & $t \ge 1$ & 3 & 1 & 15 & 2 & 20 & 1 & 96 \\
$120t + 40 +  76$ & $t \ge 1$ & 4 & 1 & 15 & 3 & 20 & 1 & 116 \\
$120t + 40 +  96$ & $t \ge 2$ & 5 & 1 & 15 & 4 & 20 & 1 & 16,\, 136,\, 256 \\
$120t + 40 + 116$ & $t \ge 2$ & 6 & 1 & 15 & 5 & 20 & 1 & 36,\, 156,\, 276
\end{tabular}
\end{center}
\end{table}

The missing values are handled as follows.
Where necessary, we give brief details by specifying the ingredients for Wilson's construction, namely
the complete graphs, the complete multipartite graphs and the group divisible designs.
Unless it is clear we also indicate how the points of the GDD are inflated and whether the GDD is augmented by an extra point.
Decompositions of the ingredients exist by
Lemmas~\ref{lem:basic multipartite direct}, \ref{lem:basic designs direct} and \ref{lem:basic designs constructed}.
For the existence of group divisible designs and mutually orthogonal Latin squares,
we usually refer the reader to \cite{Ge2007} and \cite{AbelColbournDinitz2007} respectively.



\begin{sloppypar}
There is no design of order 20 for graph $n_{13}$ by Lemma~\ref{lem:nonexistence order 20}.

Designs of order 20 for graphs $n_{3}$, $n_{6}$, $n_{8}$ and $n_{10}$ are given by Lemma~\ref{lem:basic designs direct}.

Designs of orders 40, 60, 80, 100 and 140 for all five graphs are given by
Lemmas~\ref{lem:basic designs direct} and \ref{lem:basic designs constructed}.

For order 120 for graphs $n_{3}$, $n_{6}$, $n_{8}$ and $n_{10}$, use decompositions of
$K_{20}$ and $K_{5,5,5,5,5}$ with a 5-GDD of type $4^6$ (obtained by removing a point and its incident lines from an affine plane of order 5).

For order 120 for graph $n_{13}$, use decompositions of
$K_{21}$, $K_{36}$ and $K_{21,21,21,21,36}$.

For order 180 for all five graphs, use decompositions of
$K_{45}$ and $K_{15,15,15,15}$ with a 4-GDD of type $3^4$ (obtained from a projective plane of order 3).

For order 260 for graphs $n_{3}$, $n_{6}$, $n_{8}$ and $n_{10}$, use decompositions of
$K_{20}$ and $K_{10,10,10,10}$ with a 4-GDD of type $2^{13}$, \cite{BrouwerSchrijverHanani1977}, \cite{Ge2007}.

For order 260 for graph $n_{13}$, use decompositions of
$K_{36}$, $K_{56}$, $K_{8,8,8,8,8}$ and $K_{8,8,8,8,3}$ with a 5-GDD of type $7^5$.
Inflate 4 points in one group by a factor of 3, all other points by a factor of 8.

For order 300 for graphs $n_{3}$, $n_{6}$, $n_{8}$ and $n_{10}$, use decompositions of
$K_{36}$, $K_{45}$, $K_{4^6}$ and $K_{4^6,5}$ with a $\{6,7\}$-GDD of type $11^6 7^1$
(remove 4 points from one group of a 7-GDD of type $11^7$ obtained from 5 MOLS of side 11).
Inflate the points of the reduced group by a factor of 5, all other points by a factor of 4.

For order 300 for graph $n_{13}$, use decompositions of
$K_{36}$, $K_{45}$, $K_{4^6}$, $K_{4^6,10}$ and $K_{4^6,15}$ with a $\{6,7\}$-GDD of type $11^6 3^1$
(remove 8 points from a 7-GDD of type $11^7$). 
Inflate one point in the reduced group by a factor of 15 and the other two by 10; inflate all other points by a factor of 4.


Designs of orders 21, 41, 61, 81 and 101 for all five graphs are given
by Lemmas~\ref{lem:basic designs direct} and \ref{lem:basic designs constructed}.

For order 121 for all five graphs, use decompositions of
$K_{21}$ and $K_{5,5,5,5,5}$ with a 5-GDD of type $4^6$. 

For orders 141 and 261 for all five graphs, use decompositions of
$K_{21}$ and $K_{10,10,10,10}$ with a 4-GDD of type $2^7$ or $2^{13}$, \cite{BrouwerSchrijverHanani1977}, \cite{Ge2007}.


There is no design of order 5 for any of the five graphs, by (i) of (1).

Designs of orders 25, 45, 65, 85, 105 and 125 for all five graphs are given
by Lemmas~\ref{lem:basic designs direct} and \ref{lem:basic designs constructed}.

For order 145 for all five graphs, use decompositions of
$K_{25}$ and $K_{6,6,6,6,6}$, and a 5-GDD of type $4^6$.
The GDD is augmented with an extra point.

For order 165 for all five graphs, use decompositions of
$K_{40}$, $K_{45}$, $K_{10,10,10,10}$ and $K_{10,10,10,15}$ with a $4$-GDD of type $4^4$.
Inflate one point by a factor of 15, all others by a factor of 10.

For order 245 for all five graphs, use decompositions of
$K_{40}$, $K_{45}$, $K_{10,10,10,10,10}$ and $K_{10,10,10,10,15}$ with a 5-GDD of type $4^6$.
Inflate one point by a factor of 15, all others by a factor of 10.

For order 265 for all five graphs, use decompositions of
$K_{45}$ and $K_{4^6}$ with 6-GDD of type $11^6$ obtained from 4 MOLS of side 11.
The GDD is augmented with an extra point.

For order 285 for graphs $n_{3}$, $n_{6}$, $n_{8}$ and $n_{10}$, use decompositions of
$K_{21}$, $K_{45}$, $K_{4^6}$ and $K_{4^6, 5}$ with $\{6,7\}$-GDD of type $11^6 4^1$
obtained by removing 7 points from a group of a 7-GDD of type $11^{7}$.
Inflate the points in the group of size 4 by a factor of 5, all other points by a factor of 4.
The GDD is augmented with an extra point.

For order 285 for graph $n_{13}$, use decompositions of
$K_{21}$, $K_{45}$, $K_{4^6}$ and $K_{4^6, 10}$ with a $\{6,7\}$-GDD of type $11^6 2^1$
obtained by removing 9 points from a group of a 7-GDD of type $11^{7}$.
Inflate the points in the group of size 2 by a factor of 10, all other points by a factor of 4.
The GDD is augmented with an extra point.


There is no design of order 16 for graph $n_{8}$ or $n_{13}$ by Lemma~\ref{lem:nonexistence order 16}.

We do not know if there exists a design of order 16 for graph $n_{3}$, $n_{6}$ or $n_{10}$.

Designs of orders 36, 56, 76, 96, 116 and 136 for all five graphs are given
by Lemmas~\ref{lem:basic designs direct} and \ref{lem:basic designs constructed}.

For order 156 for graphs $n_{3}$, $n_{10}$ and $n_{13}$, see Lemma~\ref{lem:basic designs direct}.

For order 156 for graphs $n_{6}$ and $n_{8}$, use decompositions of
$K_{36}$, $K_{41}$, $K_{10,10,10}$, $K_{10,10,10,10}$ and $K_{10,10,10,15}$ with a $\{3,4\}$-GDD of type $4^3 3^1$
obtained by removing a point from a 4-GDD of type $4^4$. 
Inflate one point in the reduced group by a factor of 15, all other points by a factor of 10.

For order 176 for all five graphs, use decompositions of
$K_{36}$ and $K_{5,5,5,5,5}$ with a $5$-GDD of type $7^5$.

For order 256 for graphs $n_{3}$, $n_{6}$, $n_{8}$ and $n_{10}$, use decompositions of
$K_{36}$, $K_{40}$, $K_{4^6}$ and $K_{4^6,5}$ with a $\{6,7\}$-GDD of type $9^6 8^1$
obtained by removing a point from a 7-GDD of type $9^7$ (obtained from 5 MOLS of side 9).
Inflate points in the reduced group by a factor of 5, all other points by a factor of 4.

For order 256 for graph $n_{13}$, use decompositions of
$K_{36}$, $K_{40}$, $K_{4^6}$ and $K_{4^6,10}$ with a $\{6,7\}$-GDD of type $9^6 4^1$
obtained by removing 5 points from one group of a 7-GDD of type $9^7$. 
Inflate points in the reduced group by a factor of 10, all other points by a factor of 4.

For order 276 for all five graphs, use decompositions of
$K_{56}$ and $K_{5,5,5,5,5}$ with a $5$-GDD of type $11^5$ (obtained from 3 MOLS of side 11).

For order 296 for all five graphs, use decompositions of
$K_{41}$, $K_{56}$, $K_{10,10,10,10}$ and $K_{10,10,10,15}$ with a $4$-GDD of type $4^7$ \cite{BrouwerSchrijverHanani1977}, \cite{Ge2007}.
Inflate 3 points in one group by a factor of 15, all other points by a factor of 10.
\adfQED\\
\end{sloppypar}

Thus the design spectrum for all 15 graphs with six vertices and ten edges is solved except for the
possible existence of a design of order 16 for graphs $n_3$, $n_6$ and $n_{10}$.



\appendix


\small


%
%
%
%
%
%
%
%

%
\newcommand{\ADFvfyParStart}[1]{}

\newcommand{\adfDgap}{\vskip 1.0mm}
\newcommand{\adfSgap}{\vskip 0.5mm}
\newcommand{\adfLgap}{\vskip 1.0mm}
\newcommand{\adfsplit}{\par}

\newcommand{\adfsmFJc}{3}
\newcommand{\adfsmFJf}{6}
\newcommand{\adfsmFJh}{8}
\newcommand{\adfsmFJj}{10}
\newcommand{\adfsmFJm}{13}



\section{Designs}
\label{sec:App-D}

\adfDgap
\noindent{\boldmath $ K_{20}$}~
With the point set $Z_{20}$
 the designs are generated from

\adfLgap 
$(17, 16, 13, 6, 11, 19)_{\adfsmFJc}$,


\adfLgap 
$(5, 12, 7, 13, 16, 19)_{\adfsmFJj}$,

\adfLgap \noindent by the mapping:
$x \mapsto x +  j \adfmod{19}$ for $x < 19$,
$19 \mapsto 19$,
$0 \le j < 19$.

\ADFvfyParStart{(20, ((1, 19, ((19, 1), (1, 1)))), -1)} 

\adfDgap
\noindent{\boldmath $ K_{20}$}~
With the point set $Z_{20}$
 the designs are generated from

\adfLgap 
$(19, 0, 3, 9, 7, 8)_{\adfsmFJf}$,
$(6, 0, 3, 17, 13, 2)_{\adfsmFJf}$,
$(4, 1, 0, 12, 18, 10)_{\adfsmFJf}$,\adfsplit
$(6, 12, 18, 16, 14, 19)_{\adfsmFJf}$,
$(4, 6, 11, 3, 9, 7)_{\adfsmFJf}$,
$(5, 0, 3, 15, 1, 14)_{\adfsmFJf}$,\adfsplit
$(0, 16, 5, 3, 4, 11)_{\adfsmFJf}$,
$(3, 10, 2, 12, 18, 19)_{\adfsmFJf}$,
$(1, 8, 18, 6, 7, 9)_{\adfsmFJf}$,\adfsplit
$(6, 10, 12, 5, 11, 15)_{\adfsmFJf}$,
$(8, 18, 7, 5, 12, 15)_{\adfsmFJf}$,
$(7, 10, 12, 9, 13, 17)_{\adfsmFJf}$,\adfsplit
$(8, 10, 7, 2, 14, 16)_{\adfsmFJf}$,
$(4, 8, 18, 11, 13, 17)_{\adfsmFJf}$,
$(4, 14, 11, 2, 15, 19)_{\adfsmFJf}$,\adfsplit
$(1, 11, 16, 13, 14, 17)_{\adfsmFJf}$,
$(1, 15, 5, 2, 16, 19)_{\adfsmFJf}$,
$(2, 9, 19, 13, 16, 17)_{\adfsmFJf}$,\adfsplit
$(13, 17, 9, 5, 14, 15)_{\adfsmFJf}$,


\adfLgap 
$(0, 1, 2, 6, 3, 7)_{\adfsmFJh}$,
$(2, 3, 4, 6, 5, 8)_{\adfsmFJh}$,
$(4, 5, 0, 9, 1, 10)_{\adfsmFJh}$,\adfsplit
$(18, 15, 10, 16, 7, 8)_{\adfsmFJh}$,
$(10, 8, 12, 11, 19, 17)_{\adfsmFJh}$,
$(6, 14, 16, 10, 17, 13)_{\adfsmFJh}$,\adfsplit
$(19, 9, 6, 14, 15, 18)_{\adfsmFJh}$,
$(0, 1, 8, 9, 10, 11)_{\adfsmFJh}$,
$(2, 3, 7, 10, 14, 16)_{\adfsmFJh}$,\adfsplit
$(4, 5, 6, 8, 18, 14)_{\adfsmFJh}$,
$(11, 12, 6, 15, 14, 18)_{\adfsmFJh}$,
$(0, 1, 12, 14, 13, 15)_{\adfsmFJh}$,\adfsplit
$(0, 1, 16, 17, 19, 18)_{\adfsmFJh}$,
$(2, 3, 9, 13, 19, 18)_{\adfsmFJh}$,
$(2, 3, 11, 15, 12, 17)_{\adfsmFJh}$,\adfsplit
$(4, 5, 7, 13, 12, 15)_{\adfsmFJh}$,
$(16, 17, 7, 9, 13, 12)_{\adfsmFJh}$,
$(7, 13, 8, 11, 9, 19)_{\adfsmFJh}$,\adfsplit
$(4, 5, 11, 17, 16, 19)_{\adfsmFJh}$.


\ADFvfyParStart{(20, ((19, 1, ((20, 20)))), -1)} 

\adfDgap
\noindent{\boldmath $ K_{21}$}~
With the point set $Z_{21}$
 the designs are generated from

\adfLgap 
$(0, 1, 3, 7, 13, 5)_{\adfsmFJc}$,


\adfLgap 
$(0, 1, 2, 4, 7, 10)_{\adfsmFJf}$,


\adfLgap 
$(0, 1, 2, 4, 7, 12)_{\adfsmFJh}$,


\adfLgap 
$(0, 1, 3, 7, 12, 8)_{\adfsmFJj}$,


\adfLgap 
$(0, 1, 2, 4, 7, 12)_{\adfsmFJm}$,

\adfLgap \noindent by the mapping:
$x \mapsto x +  j \adfmod{21}$,
$0 \le j < 21$.

\ADFvfyParStart{(21, ((1, 21, ((21, 1)))), -1)} 

\adfDgap
\noindent{\boldmath $ K_{25}$}~
With the point set $Z_{25}$
 the designs are generated from

\adfLgap 
$(0, 2, 5, 12, 3, 11)_{\adfsmFJc}$,
$(3, 10, 20, 1, 14, 13)_{\adfsmFJc}$,
$(13, 8, 1, 16, 22, 19)_{\adfsmFJc}$,\adfsplit
$(18, 5, 19, 14, 22, 12)_{\adfsmFJc}$,
$(19, 11, 12, 4, 17, 3)_{\adfsmFJc}$,
$(1, 4, 6, 5, 17, 14)_{\adfsmFJc}$,


\adfLgap 
$(0, 6, 1, 12, 22, 8)_{\adfsmFJf}$,
$(14, 17, 7, 9, 3, 2)_{\adfsmFJf}$,
$(5, 12, 0, 10, 18, 3)_{\adfsmFJf}$,\adfsplit
$(7, 14, 16, 6, 15, 8)_{\adfsmFJf}$,
$(4, 16, 18, 3, 13, 14)_{\adfsmFJf}$,
$(0, 11, 20, 4, 14, 16)_{\adfsmFJf}$,


\adfLgap 
$(9, 8, 5, 19, 18, 22)_{\adfsmFJh}$,
$(7, 9, 0, 10, 20, 1)_{\adfsmFJh}$,
$(2, 17, 12, 1, 18, 0)_{\adfsmFJh}$,\adfsplit
$(21, 8, 10, 15, 0, 9)_{\adfsmFJh}$,
$(1, 4, 8, 22, 12, 24)_{\adfsmFJh}$,
$(1, 23, 3, 4, 6, 11)_{\adfsmFJh}$,


\adfLgap 
$(0, 15, 3, 19, 4, 11)_{\adfsmFJj}$,
$(20, 19, 1, 11, 12, 3)_{\adfsmFJj}$,
$(19, 7, 17, 13, 18, 21)_{\adfsmFJj}$,\adfsplit
$(14, 17, 22, 1, 5, 19)_{\adfsmFJj}$,
$(9, 6, 13, 23, 11, 0)_{\adfsmFJj}$,
$(0, 2, 18, 20, 21, 22)_{\adfsmFJj}$,


\adfLgap 
$(0, 5, 12, 7, 3, 17)_{\adfsmFJm}$,
$(10, 0, 4, 9, 1, 3)_{\adfsmFJm}$,
$(5, 18, 13, 9, 1, 16)_{\adfsmFJm}$,\adfsplit
$(7, 6, 8, 13, 19, 3)_{\adfsmFJm}$,
$(7, 1, 4, 11, 21, 9)_{\adfsmFJm}$,
$(2, 5, 9, 11, 19, 13)_{\adfsmFJm}$,

\adfLgap \noindent by the mapping:
$x \mapsto x + 5 j \adfmod{25}$,
$0 \le j < 5$.

\ADFvfyParStart{(25, ((6, 5, ((25, 5)))), -1)} 

\adfDgap
\noindent{\boldmath $ K_{36}$}~
With the point set $Z_{36}$
 the designs are generated from

\adfLgap 
$(0, 20, 11, 14, 33, 6)_{\adfsmFJc}$,
$(17, 7, 28, 21, 33, 15)_{\adfsmFJc}$,
$(7, 20, 15, 19, 2, 27)_{\adfsmFJc}$,\adfsplit
$(32, 24, 34, 5, 0, 17)_{\adfsmFJc}$,
$(16, 17, 6, 23, 2, 19)_{\adfsmFJc}$,
$(23, 22, 2, 5, 14, 4)_{\adfsmFJc}$,\adfsplit
$(1, 3, 9, 10, 14, 13)_{\adfsmFJc}$,


\adfLgap 
$(0, 19, 27, 2, 31, 18)_{\adfsmFJf}$,
$(29, 0, 21, 16, 23, 6)_{\adfsmFJf}$,
$(22, 14, 13, 33, 17, 27)_{\adfsmFJf}$,\adfsplit
$(19, 12, 6, 27, 8, 20)_{\adfsmFJf}$,
$(22, 18, 25, 28, 8, 34)_{\adfsmFJf}$,
$(12, 24, 31, 13, 15, 33)_{\adfsmFJf}$,\adfsplit
$(1, 30, 17, 23, 27, 29)_{\adfsmFJf}$,


\adfLgap 
$(6, 26, 18, 7, 33, 10)_{\adfsmFJh}$,
$(8, 30, 28, 20, 33, 16)_{\adfsmFJh}$,
$(12, 28, 10, 3, 29, 21)_{\adfsmFJh}$,\adfsplit
$(9, 28, 11, 22, 34, 7)_{\adfsmFJh}$,
$(24, 23, 27, 9, 11, 1)_{\adfsmFJh}$,
$(1, 3, 4, 13, 11, 33)_{\adfsmFJh}$,\adfsplit
$(2, 6, 1, 11, 31, 16)_{\adfsmFJh}$,


\adfLgap 
$(0, 15, 16, 2, 11, 6)_{\adfsmFJj}$,
$(28, 4, 15, 31, 21, 10)_{\adfsmFJj}$,
$(19, 31, 12, 33, 2, 34)_{\adfsmFJj}$,\adfsplit
$(3, 33, 14, 25, 15, 31)_{\adfsmFJj}$,
$(23, 6, 18, 21, 14, 26)_{\adfsmFJj}$,
$(14, 34, 0, 4, 33, 1)_{\adfsmFJj}$,\adfsplit
$(1, 12, 13, 17, 4, 10)_{\adfsmFJj}$,


\adfLgap 
$(0, 11, 22, 23, 20, 18)_{\adfsmFJm}$,
$(27, 32, 21, 24, 31, 18)_{\adfsmFJm}$,
$(12, 27, 13, 33, 29, 25)_{\adfsmFJm}$,\adfsplit
$(12, 8, 1, 17, 5, 22)_{\adfsmFJm}$,
$(0, 12, 34, 26, 19, 6)_{\adfsmFJm}$,
$(14, 2, 11, 13, 21, 31)_{\adfsmFJm}$,\adfsplit
$(1, 3, 15, 10, 14, 23)_{\adfsmFJm}$,

\adfLgap \noindent by the mapping:
$x \mapsto x + 4 j \adfmod{36}$,
$0 \le j < 9$.

\ADFvfyParStart{(36, ((7, 9, ((36, 4)))), -1)} 

\adfDgap
\noindent{\boldmath $ K_{40}$}~
With the point set $Z_{40}$
 the designs are generated from

\adfLgap 
$(0, 24, 35, 3, 29, 39)_{\adfsmFJc}$,
$(0, 8, 22, 9, 20, 16)_{\adfsmFJc}$,


\adfLgap 
$(2, 0, 38, 27, 18, 39)_{\adfsmFJj}$,
$(0, 4, 17, 33, 9, 7)_{\adfsmFJj}$,

\adfLgap \noindent by the mapping:
$x \mapsto x +  j \adfmod{39}$ for $x < 39$,
$39 \mapsto 39$,
$0 \le j < 39$.

\ADFvfyParStart{(40, ((2, 39, ((39, 1), (1, 1)))), -1)} 

\adfDgap
\noindent{\boldmath $ K_{40}$}~
With the point set $Z_{40}$
 the designs are generated from

\adfLgap 
$(0, 17, 39, 29, 6, 10)_{\adfsmFJf}$,
$(4, 8, 17, 26, 9, 14)_{\adfsmFJf}$,
$(35, 16, 10, 30, 9, 33)_{\adfsmFJf}$,\adfsplit
$(8, 24, 3, 15, 16, 27)_{\adfsmFJf}$,
$(27, 2, 10, 6, 16, 25)_{\adfsmFJf}$,
$(13, 26, 22, 11, 25, 31)_{\adfsmFJf}$,


\adfLgap 
$(18, 14, 39, 10, 37, 12)_{\adfsmFJh}$,
$(15, 35, 10, 22, 29, 6)_{\adfsmFJh}$,
$(2, 20, 13, 10, 19, 9)_{\adfsmFJh}$,\adfsplit
$(10, 34, 21, 31, 24, 19)_{\adfsmFJh}$,
$(31, 14, 27, 29, 9, 26)_{\adfsmFJh}$,
$(36, 2, 20, 24, 21, 32)_{\adfsmFJh}$,


\adfLgap 
$(0, 39, 29, 14, 34, 36)_{\adfsmFJm}$,
$(18, 14, 7, 30, 10, 38)_{\adfsmFJm}$,
$(38, 22, 8, 23, 1, 5)_{\adfsmFJm}$,\adfsplit
$(5, 32, 34, 19, 33, 25)_{\adfsmFJm}$,
$(32, 27, 30, 6, 4, 10)_{\adfsmFJm}$,
$(0, 10, 30, 22, 24, 37)_{\adfsmFJm}$,

\adfLgap \noindent by the mapping:
$x \mapsto x + 3 j \adfmod{39}$ for $x < 39$,
$39 \mapsto 39$,
$0 \le j < 13$.

\ADFvfyParStart{(40, ((6, 13, ((39, 3), (1, 1)))), -1)} 

\adfDgap
\noindent{\boldmath $ K_{41}$}~
With the point set $Z_{41}$
 the designs are generated from

\adfLgap 
$(0, 1, 3, 7, 13, 15)_{\adfsmFJc}$,
$(0, 5, 14, 22, 25, 18)_{\adfsmFJc}$,


\adfLgap 
$(0, 1, 2, 4, 7, 10)_{\adfsmFJf}$,
$(0, 11, 3, 24, 25, 29)_{\adfsmFJf}$,


\adfLgap 
$(0, 1, 2, 4, 7, 15)_{\adfsmFJh}$,
$(0, 4, 13, 16, 23, 24)_{\adfsmFJh}$,


\adfLgap 
$(0, 1, 3, 7, 16, 11)_{\adfsmFJj}$,
$(0, 8, 20, 25, 39, 18)_{\adfsmFJj}$,


\adfLgap 
$(0, 1, 2, 4, 7, 17)_{\adfsmFJm}$,
$(0, 9, 11, 19, 29, 27)_{\adfsmFJm}$,

\adfLgap \noindent by the mapping:
$x \mapsto x +  j \adfmod{41}$,
$0 \le j < 41$.

\ADFvfyParStart{(41, ((2, 41, ((41, 1)))), -1)} 

\adfDgap
\noindent{\boldmath $ K_{45}$}~
With the point set $Z_{45}$
 the designs are generated from

\adfLgap 
$(37, 44, 43, 8, 10, 4)_{\adfsmFJc}$,
$(27, 13, 35, 23, 6, 3)_{\adfsmFJc}$,
$(1, 36, 5, 26, 28, 4)_{\adfsmFJc}$,\adfsplit
$(5, 17, 7, 14, 12, 23)_{\adfsmFJc}$,
$(25, 38, 33, 9, 19, 24)_{\adfsmFJc}$,
$(35, 28, 30, 34, 22, 9)_{\adfsmFJc}$,\adfsplit
$(30, 41, 16, 42, 39, 8)_{\adfsmFJc}$,
$(38, 20, 18, 8, 35, 10)_{\adfsmFJc}$,
$(3, 0, 4, 28, 31, 36)_{\adfsmFJc}$,


\adfLgap 
$(44, 42, 20, 29, 16, 19)_{\adfsmFJf}$,
$(38, 28, 18, 2, 23, 9)_{\adfsmFJf}$,
$(16, 17, 5, 43, 38, 22)_{\adfsmFJf}$,\adfsplit
$(33, 41, 22, 34, 11, 20)_{\adfsmFJf}$,
$(36, 6, 22, 23, 25, 26)_{\adfsmFJf}$,
$(39, 14, 13, 11, 16, 28)_{\adfsmFJf}$,\adfsplit
$(38, 7, 28, 0, 31, 3)_{\adfsmFJf}$,
$(29, 1, 43, 31, 33, 35)_{\adfsmFJf}$,
$(0, 20, 17, 12, 35, 37)_{\adfsmFJf}$,


\adfLgap 
$(44, 11, 28, 31, 26, 41)_{\adfsmFJh}$,
$(26, 1, 33, 8, 18, 19)_{\adfsmFJh}$,
$(0, 20, 6, 39, 12, 41)_{\adfsmFJh}$,\adfsplit
$(34, 30, 14, 9, 31, 39)_{\adfsmFJh}$,
$(22, 29, 12, 35, 40, 25)_{\adfsmFJh}$,
$(9, 29, 31, 0, 37, 4)_{\adfsmFJh}$,\adfsplit
$(26, 33, 38, 13, 15, 22)_{\adfsmFJh}$,
$(15, 12, 3, 36, 11, 43)_{\adfsmFJh}$,
$(0, 12, 13, 15, 14, 34)_{\adfsmFJh}$,


\adfLgap 
$(9, 16, 44, 39, 38, 22)_{\adfsmFJj}$,
$(26, 9, 28, 8, 35, 25)_{\adfsmFJj}$,
$(25, 28, 33, 30, 18, 43)_{\adfsmFJj}$,\adfsplit
$(14, 39, 5, 25, 3, 11)_{\adfsmFJj}$,
$(35, 23, 24, 19, 6, 20)_{\adfsmFJj}$,
$(7, 31, 42, 2, 21, 32)_{\adfsmFJj}$,\adfsplit
$(39, 18, 2, 8, 38, 36)_{\adfsmFJj}$,
$(27, 25, 21, 36, 9, 33)_{\adfsmFJj}$,
$(0, 9, 32, 40, 2, 16)_{\adfsmFJj}$,

\adfLgap \noindent by the mapping:
$x \mapsto x + 4 j \adfmod{44}$ for $x < 44$,
$44 \mapsto 44$,
$0 \le j < 11$.

\ADFvfyParStart{(45, ((9, 11, ((44, 4), (1, 1)))), -1)} 

\adfDgap
\noindent{\boldmath $ K_{45}$}~
With the point set $Z_{45}$
 the design is generated from

\adfLgap 
$(0, 15, 21, 26, 25, 12)_{\adfsmFJm}$,
$(42, 23, 18, 34, 17, 20)_{\adfsmFJm}$,
$(37, 36, 34, 8, 7, 35)_{\adfsmFJm}$,\adfsplit
$(37, 11, 32, 5, 23, 39)_{\adfsmFJm}$,
$(5, 12, 22, 8, 34, 9)_{\adfsmFJm}$,
$(5, 13, 0, 2, 33, 36)_{\adfsmFJm}$,\adfsplit
$(0, 24, 38, 19, 16, 1)_{\adfsmFJm}$,
$(37, 27, 16, 0, 31, 9)_{\adfsmFJm}$,
$(1, 3, 4, 8, 6, 21)_{\adfsmFJm}$,\adfsplit
$(3, 19, 29, 25, 33, 38)_{\adfsmFJm}$,
$(4, 16, 19, 26, 43, 39)_{\adfsmFJm}$,

\adfLgap \noindent by the mapping:
$x \mapsto x + 5 j \adfmod{45}$,
$0 \le j < 9$.

\ADFvfyParStart{(45, ((11, 9, ((45, 5)))), -1)} 

\adfDgap
\noindent{\boldmath $ K_{56}$}~
With the point set $Z_{56}$
 the designs are generated from

\adfLgap 
$(55, 23, 49, 15, 6, 52)_{\adfsmFJc}$,
$(4, 11, 19, 45, 21, 3)_{\adfsmFJc}$,
$(52, 31, 34, 23, 53, 11)_{\adfsmFJc}$,\adfsplit
$(38, 51, 32, 29, 0, 6)_{\adfsmFJc}$,
$(47, 31, 46, 40, 5, 21)_{\adfsmFJc}$,
$(26, 17, 28, 38, 25, 22)_{\adfsmFJc}$,\adfsplit
$(54, 41, 23, 9, 38, 29)_{\adfsmFJc}$,
$(21, 32, 1, 28, 26, 39)_{\adfsmFJc}$,
$(23, 45, 7, 42, 5, 50)_{\adfsmFJc}$,\adfsplit
$(11, 39, 35, 0, 30, 50)_{\adfsmFJc}$,
$(48, 18, 25, 13, 14, 50)_{\adfsmFJc}$,
$(0, 13, 27, 19, 22, 10)_{\adfsmFJc}$,\adfsplit
$(2, 4, 17, 24, 27, 29)_{\adfsmFJc}$,
$(2, 14, 20, 19, 45, 33)_{\adfsmFJc}$,


\adfLgap 
$(55, 51, 4, 48, 37, 17)_{\adfsmFJj}$,
$(30, 11, 37, 7, 54, 55)_{\adfsmFJj}$,
$(34, 20, 33, 35, 13, 11)_{\adfsmFJj}$,\adfsplit
$(44, 16, 15, 46, 18, 23)_{\adfsmFJj}$,
$(6, 26, 27, 10, 16, 22)_{\adfsmFJj}$,
$(42, 40, 23, 32, 27, 29)_{\adfsmFJj}$,\adfsplit
$(45, 10, 40, 34, 49, 27)_{\adfsmFJj}$,
$(46, 31, 22, 23, 49, 27)_{\adfsmFJj}$,
$(1, 35, 7, 14, 23, 39)_{\adfsmFJj}$,\adfsplit
$(30, 1, 38, 34, 15, 41)_{\adfsmFJj}$,
$(6, 18, 13, 28, 42, 11)_{\adfsmFJj}$,
$(39, 46, 8, 33, 49, 29)_{\adfsmFJj}$,\adfsplit
$(14, 34, 2, 39, 37, 33)_{\adfsmFJj}$,
$(0, 10, 22, 28, 7, 46)_{\adfsmFJj}$,

\adfLgap \noindent by the mapping:
$x \mapsto x + 5 j \adfmod{55}$ for $x < 55$,
$55 \mapsto 55$,
$0 \le j < 11$.

\ADFvfyParStart{(56, ((14, 11, ((55, 5), (1, 1)))), -1)} 

\adfDgap
\noindent{\boldmath $ K_{56}$}~
With the point set $Z_{56}$
 the designs are generated from

\adfLgap 
$(18, 23, 47, 4, 49, 38)_{\adfsmFJf}$,
$(20, 25, 49, 6, 51, 40)_{\adfsmFJf}$,
$(22, 27, 51, 8, 53, 42)_{\adfsmFJf}$,\adfsplit
$(24, 29, 53, 10, 55, 44)_{\adfsmFJf}$,
$(50, 7, 6, 51, 1, 53)_{\adfsmFJf}$,
$(52, 9, 8, 53, 3, 55)_{\adfsmFJf}$,\adfsplit
$(54, 11, 10, 55, 5, 1)_{\adfsmFJf}$,
$(0, 13, 12, 1, 7, 3)_{\adfsmFJf}$,
$(23, 45, 29, 6, 8, 52)_{\adfsmFJf}$,\adfsplit
$(25, 47, 31, 8, 10, 54)_{\adfsmFJf}$,
$(27, 49, 33, 10, 12, 0)_{\adfsmFJf}$,
$(29, 51, 35, 12, 14, 2)_{\adfsmFJf}$,\adfsplit
$(12, 42, 40, 44, 24, 34)_{\adfsmFJf}$,
$(14, 44, 42, 46, 26, 36)_{\adfsmFJf}$,
$(16, 46, 44, 48, 28, 38)_{\adfsmFJf}$,\adfsplit
$(18, 48, 46, 50, 30, 40)_{\adfsmFJf}$,
$(53, 29, 46, 15, 0, 45)_{\adfsmFJf}$,
$(55, 31, 48, 17, 2, 47)_{\adfsmFJf}$,\adfsplit
$(15, 35, 47, 43, 36, 51)_{\adfsmFJf}$,
$(17, 41, 34, 3, 6, 44)_{\adfsmFJf}$,
$(1, 37, 13, 2, 9, 17)_{\adfsmFJf}$,\adfsplit
$(3, 27, 20, 30, 45, 48)_{\adfsmFJf}$,


\adfLgap 
$(19, 37, 50, 13, 54, 0)_{\adfsmFJh}$,
$(21, 39, 52, 15, 0, 2)_{\adfsmFJh}$,
$(23, 41, 54, 17, 2, 4)_{\adfsmFJh}$,\adfsplit
$(25, 43, 0, 19, 4, 6)_{\adfsmFJh}$,
$(28, 1, 54, 50, 8, 42)_{\adfsmFJh}$,
$(30, 3, 0, 52, 10, 44)_{\adfsmFJh}$,\adfsplit
$(32, 5, 2, 54, 12, 46)_{\adfsmFJh}$,
$(34, 7, 4, 0, 14, 48)_{\adfsmFJh}$,
$(15, 14, 54, 25, 45, 16)_{\adfsmFJh}$,\adfsplit
$(17, 16, 0, 27, 47, 18)_{\adfsmFJh}$,
$(19, 18, 2, 29, 49, 20)_{\adfsmFJh}$,
$(21, 20, 4, 31, 51, 22)_{\adfsmFJh}$,\adfsplit
$(48, 36, 4, 25, 33, 54)_{\adfsmFJh}$,
$(50, 38, 6, 27, 35, 0)_{\adfsmFJh}$,
$(52, 40, 8, 29, 37, 2)_{\adfsmFJh}$,\adfsplit
$(54, 42, 10, 31, 39, 4)_{\adfsmFJh}$,
$(37, 11, 49, 33, 51, 32)_{\adfsmFJh}$,
$(39, 13, 51, 35, 53, 34)_{\adfsmFJh}$,\adfsplit
$(41, 15, 53, 37, 55, 36)_{\adfsmFJh}$,
$(43, 17, 55, 39, 1, 38)_{\adfsmFJh}$,
$(1, 29, 6, 9, 34, 37)_{\adfsmFJh}$,\adfsplit
$(3, 31, 8, 11, 36, 39)_{\adfsmFJh}$,


\adfLgap 
$(20, 12, 44, 14, 17, 51)_{\adfsmFJm}$,
$(25, 27, 31, 44, 41, 3)_{\adfsmFJm}$,
$(3, 10, 50, 38, 11, 18)_{\adfsmFJm}$,\adfsplit
$(54, 1, 46, 21, 50, 25)_{\adfsmFJm}$,
$(19, 6, 1, 3, 0, 9)_{\adfsmFJm}$,
$(11, 29, 8, 21, 23, 32)_{\adfsmFJm}$,\adfsplit
$(24, 45, 5, 14, 11, 4)_{\adfsmFJm}$,
$(5, 31, 17, 36, 54, 42)_{\adfsmFJm}$,
$(32, 25, 34, 23, 26, 37)_{\adfsmFJm}$,\adfsplit
$(29, 25, 36, 30, 49, 13)_{\adfsmFJm}$,
$(1, 23, 31, 13, 34, 18)_{\adfsmFJm}$,
$(15, 46, 54, 37, 10, 13)_{\adfsmFJm}$,\adfsplit
$(24, 41, 31, 54, 55, 38)_{\adfsmFJm}$,
$(55, 11, 20, 15, 37, 16)_{\adfsmFJm}$,
$(40, 22, 36, 10, 52, 50)_{\adfsmFJm}$,\adfsplit
$(35, 44, 40, 6, 9, 11)_{\adfsmFJm}$,
$(40, 38, 32, 6, 4, 17)_{\adfsmFJm}$,
$(2, 16, 28, 45, 39, 17)_{\adfsmFJm}$,\adfsplit
$(26, 36, 8, 55, 48, 17)_{\adfsmFJm}$,
$(4, 42, 27, 13, 7, 2)_{\adfsmFJm}$,
$(3, 4, 5, 45, 26, 48)_{\adfsmFJm}$,\adfsplit
$(54, 43, 55, 39, 7, 12)_{\adfsmFJm}$,

\adfLgap \noindent by the mapping:
$x \mapsto x + 8 j \adfmod{56}$,
$0 \le j < 7$.

\ADFvfyParStart{(56, ((22, 7, ((56, 8)))), -1)} 

\adfDgap
\noindent{\boldmath $ K_{60}$}~
With the point set $Z_{60}$
 the designs are generated from

\adfLgap 
$(22, 15, 10, 2, 43, 59)_{\adfsmFJc}$,
$(0, 1, 3, 18, 35, 9)_{\adfsmFJc}$,
$(0, 6, 29, 10, 43, 11)_{\adfsmFJc}$,


\adfLgap 
$(4, 13, 50, 19, 53, 59)_{\adfsmFJj}$,
$(45, 4, 16, 40, 42, 0)_{\adfsmFJj}$,
$(0, 1, 17, 49, 56, 8)_{\adfsmFJj}$,

\adfLgap \noindent by the mapping:
$x \mapsto x +  j \adfmod{59}$ for $x < 59$,
$59 \mapsto 59$,
$0 \le j < 59$.

\ADFvfyParStart{(60, ((3, 59, ((59, 1), (1, 1)))), -1)} 

\adfDgap
\noindent{\boldmath $ K_{61}$}~
With the point set $Z_{61}$
 the designs are generated from

\adfLgap 
$(0, 47, 32, 34, 8, 50)_{\adfsmFJc}$,
$(29, 60, 50, 17, 24, 10)_{\adfsmFJc}$,
$(22, 5, 6, 2, 60, 13)_{\adfsmFJc}$,


\adfLgap 
$(0, 16, 26, 31, 50, 3)_{\adfsmFJj}$,
$(45, 46, 53, 5, 14, 1)_{\adfsmFJj}$,
$(0, 4, 6, 18, 29, 28)_{\adfsmFJj}$,


\adfLgap 
$(0, 29, 11, 14, 13, 37)_{\adfsmFJm}$,
$(50, 9, 38, 6, 8, 45)_{\adfsmFJm}$,
$(0, 6, 28, 27, 10, 36)_{\adfsmFJm}$,

\adfLgap \noindent by the mapping:
$x \mapsto x +  j \adfmod{61}$,
$0 \le j < 61$.

\ADFvfyParStart{(61, ((3, 61, ((61, 1)))), -1)} 

\adfDgap
\noindent{\boldmath $ K_{65}$}~
With the point set $Z_{65}$
 the designs are generated from

\adfLgap 
$(30, 12, 60, 49, 53, 27)_{\adfsmFJc}$,
$(31, 13, 61, 50, 54, 28)_{\adfsmFJc}$,
$(32, 14, 62, 51, 55, 29)_{\adfsmFJc}$,\adfsplit
$(33, 15, 63, 52, 56, 30)_{\adfsmFJc}$,
$(34, 16, 64, 53, 57, 31)_{\adfsmFJc}$,
$(55, 5, 10, 54, 43, 46)_{\adfsmFJc}$,\adfsplit
$(56, 6, 11, 55, 44, 47)_{\adfsmFJc}$,
$(57, 7, 12, 56, 45, 48)_{\adfsmFJc}$,
$(58, 8, 13, 57, 46, 49)_{\adfsmFJc}$,\adfsplit
$(59, 9, 14, 58, 47, 50)_{\adfsmFJc}$,
$(30, 32, 36, 61, 22, 8)_{\adfsmFJc}$,
$(31, 33, 37, 62, 23, 9)_{\adfsmFJc}$,\adfsplit
$(32, 34, 38, 63, 24, 10)_{\adfsmFJc}$,
$(48, 22, 35, 61, 9, 26)_{\adfsmFJc}$,
$(1, 5, 64, 30, 56, 14)_{\adfsmFJc}$,\adfsplit
$(4, 0, 63, 29, 55, 47)_{\adfsmFJc}$,


\adfLgap 
$(61, 56, 1, 17, 25, 49)_{\adfsmFJf}$,
$(41, 6, 43, 51, 44, 42)_{\adfsmFJf}$,
$(19, 0, 15, 55, 9, 11)_{\adfsmFJf}$,\adfsplit
$(46, 0, 61, 28, 58, 23)_{\adfsmFJf}$,
$(51, 58, 52, 17, 3, 47)_{\adfsmFJf}$,
$(4, 32, 64, 38, 20, 17)_{\adfsmFJf}$,\adfsplit
$(44, 32, 45, 53, 42, 28)_{\adfsmFJf}$,
$(54, 47, 50, 18, 55, 0)_{\adfsmFJf}$,
$(10, 26, 11, 13, 17, 44)_{\adfsmFJf}$,\adfsplit
$(24, 39, 45, 0, 18, 43)_{\adfsmFJf}$,
$(0, 2, 39, 35, 26, 22)_{\adfsmFJf}$,
$(29, 52, 21, 15, 49, 43)_{\adfsmFJf}$,\adfsplit
$(59, 19, 46, 27, 24, 61)_{\adfsmFJf}$,
$(28, 4, 8, 58, 47, 42)_{\adfsmFJf}$,
$(52, 12, 26, 25, 35, 1)_{\adfsmFJf}$,\adfsplit
$(18, 43, 51, 23, 0, 30)_{\adfsmFJf}$,


\adfLgap 
$(36, 25, 49, 22, 54, 16)_{\adfsmFJh}$,
$(1, 10, 23, 61, 51, 55)_{\adfsmFJh}$,
$(39, 52, 0, 61, 12, 31)_{\adfsmFJh}$,\adfsplit
$(25, 43, 62, 41, 7, 8)_{\adfsmFJh}$,
$(34, 43, 61, 19, 64, 35)_{\adfsmFJh}$,
$(45, 37, 43, 15, 33, 11)_{\adfsmFJh}$,\adfsplit
$(2, 14, 0, 62, 34, 21)_{\adfsmFJh}$,
$(49, 40, 63, 57, 21, 61)_{\adfsmFJh}$,
$(35, 42, 50, 62, 3, 39)_{\adfsmFJh}$,\adfsplit
$(13, 64, 27, 33, 57, 39)_{\adfsmFJh}$,
$(46, 52, 36, 45, 43, 20)_{\adfsmFJh}$,
$(32, 8, 44, 21, 61, 19)_{\adfsmFJh}$,\adfsplit
$(51, 40, 17, 11, 19, 3)_{\adfsmFJh}$,
$(14, 63, 60, 25, 5, 20)_{\adfsmFJh}$,
$(22, 48, 21, 53, 44, 37)_{\adfsmFJh}$,\adfsplit
$(48, 54, 7, 64, 8, 63)_{\adfsmFJh}$,


\adfLgap 
$(54, 49, 24, 13, 47, 28)_{\adfsmFJj}$,
$(55, 50, 25, 14, 48, 29)_{\adfsmFJj}$,
$(56, 51, 26, 15, 49, 30)_{\adfsmFJj}$,\adfsplit
$(57, 52, 27, 16, 50, 31)_{\adfsmFJj}$,
$(58, 53, 28, 17, 51, 32)_{\adfsmFJj}$,
$(41, 55, 47, 54, 10, 51)_{\adfsmFJj}$,\adfsplit
$(42, 56, 48, 55, 11, 52)_{\adfsmFJj}$,
$(43, 57, 49, 56, 12, 53)_{\adfsmFJj}$,
$(44, 58, 50, 57, 13, 54)_{\adfsmFJj}$,\adfsplit
$(45, 59, 51, 58, 14, 55)_{\adfsmFJj}$,
$(41, 37, 22, 19, 10, 24)_{\adfsmFJj}$,
$(42, 38, 23, 20, 11, 25)_{\adfsmFJj}$,\adfsplit
$(43, 39, 24, 21, 12, 26)_{\adfsmFJj}$,
$(49, 58, 11, 61, 15, 17)_{\adfsmFJj}$,
$(3, 12, 15, 30, 34, 19)_{\adfsmFJj}$,\adfsplit
$(3, 20, 36, 52, 4, 35)_{\adfsmFJj}$,


\adfLgap 
$(28, 50, 0, 2, 22, 35)_{\adfsmFJm}$,
$(4, 22, 19, 31, 15, 11)_{\adfsmFJm}$,
$(29, 23, 58, 54, 46, 60)_{\adfsmFJm}$,\adfsplit
$(31, 64, 61, 37, 27, 40)_{\adfsmFJm}$,
$(34, 11, 37, 18, 35, 47)_{\adfsmFJm}$,
$(28, 62, 25, 53, 17, 20)_{\adfsmFJm}$,\adfsplit
$(31, 44, 46, 18, 20, 56)_{\adfsmFJm}$,
$(19, 9, 41, 54, 33, 52)_{\adfsmFJm}$,
$(55, 34, 3, 43, 13, 6)_{\adfsmFJm}$,\adfsplit
$(32, 11, 45, 57, 16, 31)_{\adfsmFJm}$,
$(42, 12, 58, 34, 63, 57)_{\adfsmFJm}$,
$(3, 27, 64, 22, 50, 46)_{\adfsmFJm}$,\adfsplit
$(8, 43, 9, 41, 28, 35)_{\adfsmFJm}$,
$(39, 22, 36, 51, 34, 55)_{\adfsmFJm}$,
$(40, 20, 55, 49, 2, 30)_{\adfsmFJm}$,\adfsplit
$(38, 42, 55, 49, 56, 21)_{\adfsmFJm}$,

\adfLgap \noindent by the mapping:
$x \mapsto x + 5 j \adfmod{65}$,
$0 \le j < 13$.

\ADFvfyParStart{(65, ((16, 13, ((65, 5)))), -1)} 

\adfDgap
\noindent{\boldmath $ K_{76}$}~
With the point set $Z_{76}$
 the designs are generated from

\adfLgap 
$(0, 33, 42, 5, 44, 41)_{\adfsmFJc}$,
$(64, 6, 46, 27, 36, 25)_{\adfsmFJc}$,
$(52, 12, 25, 19, 66, 1)_{\adfsmFJc}$,\adfsplit
$(48, 42, 49, 67, 19, 28)_{\adfsmFJc}$,
$(13, 38, 65, 71, 42, 29)_{\adfsmFJc}$,
$(10, 71, 56, 42, 58, 27)_{\adfsmFJc}$,\adfsplit
$(34, 1, 14, 23, 8, 58)_{\adfsmFJc}$,
$(74, 57, 17, 13, 43, 14)_{\adfsmFJc}$,
$(50, 53, 43, 42, 55, 49)_{\adfsmFJc}$,\adfsplit
$(1, 46, 11, 43, 56, 6)_{\adfsmFJc}$,
$(17, 29, 0, 9, 67, 31)_{\adfsmFJc}$,
$(6, 18, 60, 55, 44, 61)_{\adfsmFJc}$,\adfsplit
$(17, 15, 32, 40, 36, 20)_{\adfsmFJc}$,
$(0, 3, 45, 23, 75, 12)_{\adfsmFJc}$,
$(0, 27, 35, 11, 63, 24)_{\adfsmFJc}$,


\adfLgap 
$(59, 1, 21, 30, 73, 9)_{\adfsmFJf}$,
$(26, 23, 74, 44, 66, 39)_{\adfsmFJf}$,
$(66, 47, 71, 29, 51, 44)_{\adfsmFJf}$,\adfsplit
$(2, 54, 68, 3, 73, 33)_{\adfsmFJf}$,
$(48, 23, 24, 73, 55, 47)_{\adfsmFJf}$,
$(52, 45, 44, 42, 73, 15)_{\adfsmFJf}$,\adfsplit
$(45, 9, 51, 39, 23, 2)_{\adfsmFJf}$,
$(42, 74, 29, 48, 35, 51)_{\adfsmFJf}$,
$(2, 28, 32, 47, 65, 66)_{\adfsmFJf}$,\adfsplit
$(49, 33, 46, 36, 74, 50)_{\adfsmFJf}$,
$(29, 19, 38, 60, 61, 27)_{\adfsmFJf}$,
$(37, 39, 4, 75, 48, 52)_{\adfsmFJf}$,\adfsplit
$(68, 70, 65, 56, 64, 9)_{\adfsmFJf}$,
$(13, 2, 6, 23, 62, 36)_{\adfsmFJf}$,
$(12, 28, 66, 48, 4, 71)_{\adfsmFJf}$,


\adfLgap 
$(71, 34, 0, 73, 42, 21)_{\adfsmFJh}$,
$(62, 28, 75, 60, 48, 65)_{\adfsmFJh}$,
$(51, 16, 17, 61, 29, 12)_{\adfsmFJh}$,\adfsplit
$(33, 20, 43, 12, 75, 53)_{\adfsmFJh}$,
$(62, 22, 71, 11, 53, 19)_{\adfsmFJh}$,
$(20, 65, 46, 3, 56, 36)_{\adfsmFJh}$,\adfsplit
$(71, 6, 41, 66, 47, 17)_{\adfsmFJh}$,
$(14, 21, 15, 7, 54, 67)_{\adfsmFJh}$,
$(61, 8, 66, 2, 38, 6)_{\adfsmFJh}$,\adfsplit
$(22, 32, 43, 2, 7, 54)_{\adfsmFJh}$,
$(69, 22, 21, 12, 37, 73)_{\adfsmFJh}$,
$(51, 32, 60, 49, 47, 31)_{\adfsmFJh}$,\adfsplit
$(57, 72, 60, 21, 49, 50)_{\adfsmFJh}$,
$(19, 66, 16, 60, 47, 2)_{\adfsmFJh}$,
$(51, 20, 39, 44, 13, 6)_{\adfsmFJh}$,


\adfLgap 
$(70, 1, 12, 52, 6, 32)_{\adfsmFJj}$,
$(73, 15, 28, 43, 32, 52)_{\adfsmFJj}$,
$(47, 15, 46, 50, 61, 71)_{\adfsmFJj}$,\adfsplit
$(27, 1, 25, 61, 35, 4)_{\adfsmFJj}$,
$(36, 29, 56, 73, 32, 44)_{\adfsmFJj}$,
$(72, 5, 9, 73, 67, 75)_{\adfsmFJj}$,\adfsplit
$(34, 44, 70, 73, 19, 68)_{\adfsmFJj}$,
$(46, 18, 38, 73, 50, 60)_{\adfsmFJj}$,
$(10, 43, 47, 70, 55, 41)_{\adfsmFJj}$,\adfsplit
$(67, 50, 72, 74, 40, 69)_{\adfsmFJj}$,
$(5, 20, 25, 34, 53, 18)_{\adfsmFJj}$,
$(20, 59, 70, 75, 69, 63)_{\adfsmFJj}$,\adfsplit
$(40, 56, 62, 75, 37, 52)_{\adfsmFJj}$,
$(30, 11, 51, 73, 20, 13)_{\adfsmFJj}$,
$(32, 2, 4, 31, 11, 3)_{\adfsmFJj}$,


\adfLgap 
$(19, 73, 53, 33, 6, 20)_{\adfsmFJm}$,
$(27, 55, 72, 4, 21, 48)_{\adfsmFJm}$,
$(74, 44, 46, 39, 32, 41)_{\adfsmFJm}$,\adfsplit
$(31, 16, 67, 7, 35, 13)_{\adfsmFJm}$,
$(36, 1, 74, 53, 33, 16)_{\adfsmFJm}$,
$(3, 18, 21, 74, 11, 0)_{\adfsmFJm}$,\adfsplit
$(67, 53, 32, 7, 34, 40)_{\adfsmFJm}$,
$(30, 36, 46, 41, 22, 40)_{\adfsmFJm}$,
$(61, 6, 1, 31, 2, 54)_{\adfsmFJm}$,\adfsplit
$(36, 7, 54, 62, 5, 17)_{\adfsmFJm}$,
$(8, 30, 67, 42, 56, 17)_{\adfsmFJm}$,
$(69, 41, 8, 61, 54, 57)_{\adfsmFJm}$,\adfsplit
$(37, 60, 48, 47, 44, 8)_{\adfsmFJm}$,
$(39, 50, 10, 51, 19, 13)_{\adfsmFJm}$,
$(26, 43, 23, 70, 62, 16)_{\adfsmFJm}$,

\adfLgap \noindent by the mapping:
$x \mapsto x + 4 j \adfmod{76}$,
$0 \le j < 19$.

\ADFvfyParStart{(76, ((15, 19, ((76, 4)))), -1)} 

\adfDgap
\noindent{\boldmath $ K_{80}$}~
With the point set $Z_{80}$
 the designs are generated from

\adfLgap 
$(43, 74, 10, 22, 69, 59)_{\adfsmFJc}$,
$(60, 3, 74, 73, 63, 58)_{\adfsmFJc}$,
$(77, 2, 36, 42, 26, 21)_{\adfsmFJc}$,\adfsplit
$(18, 11, 43, 60, 61, 79)_{\adfsmFJc}$,


\adfLgap 
$(22, 31, 55, 14, 16, 59)_{\adfsmFJj}$,
$(22, 67, 53, 17, 78, 38)_{\adfsmFJj}$,
$(52, 72, 25, 46, 76, 51)_{\adfsmFJj}$,\adfsplit
$(56, 0, 12, 69, 72, 79)_{\adfsmFJj}$,

\adfLgap \noindent by the mapping:
$x \mapsto x +  j \adfmod{79}$ for $x < 79$,
$79 \mapsto 79$,
$0 \le j < 79$.

\ADFvfyParStart{(80, ((4, 79, ((79, 1), (1, 1)))), -1)} 

\adfDgap
\noindent{\boldmath $ K_{85}$}~
With the point set $Z_{85}$
 the designs are generated from

\adfLgap 
$(32, 64, 67, 20, 63, 25)_{\adfsmFJc}$,
$(33, 65, 68, 21, 64, 26)_{\adfsmFJc}$,
$(34, 66, 69, 22, 65, 27)_{\adfsmFJc}$,\adfsplit
$(35, 67, 70, 23, 66, 28)_{\adfsmFJc}$,
$(36, 68, 71, 24, 67, 29)_{\adfsmFJc}$,
$(8, 17, 72, 33, 66, 6)_{\adfsmFJc}$,\adfsplit
$(9, 18, 73, 34, 67, 7)_{\adfsmFJc}$,
$(10, 19, 74, 35, 68, 8)_{\adfsmFJc}$,
$(11, 20, 75, 36, 69, 9)_{\adfsmFJc}$,\adfsplit
$(12, 21, 76, 37, 70, 10)_{\adfsmFJc}$,
$(73, 81, 55, 44, 10, 7)_{\adfsmFJc}$,
$(74, 82, 56, 45, 11, 8)_{\adfsmFJc}$,\adfsplit
$(75, 83, 57, 46, 12, 9)_{\adfsmFJc}$,
$(76, 84, 58, 47, 13, 10)_{\adfsmFJc}$,
$(77, 0, 59, 48, 14, 11)_{\adfsmFJc}$,\adfsplit
$(5, 25, 67, 77, 10, 29)_{\adfsmFJc}$,
$(6, 26, 68, 78, 11, 30)_{\adfsmFJc}$,
$(7, 27, 69, 79, 12, 31)_{\adfsmFJc}$,\adfsplit
$(8, 28, 70, 80, 13, 32)_{\adfsmFJc}$,
$(3, 20, 54, 37, 71, 64)_{\adfsmFJc}$,
$(1, 14, 34, 29, 76, 52)_{\adfsmFJc}$,


\adfLgap 
$(62, 48, 19, 35, 50, 22)_{\adfsmFJj}$,
$(63, 49, 20, 36, 51, 23)_{\adfsmFJj}$,
$(64, 50, 21, 37, 52, 24)_{\adfsmFJj}$,\adfsplit
$(65, 51, 22, 38, 53, 25)_{\adfsmFJj}$,
$(66, 52, 23, 39, 54, 26)_{\adfsmFJj}$,
$(9, 37, 56, 48, 47, 72)_{\adfsmFJj}$,\adfsplit
$(10, 38, 57, 49, 48, 73)_{\adfsmFJj}$,
$(11, 39, 58, 50, 49, 74)_{\adfsmFJj}$,
$(12, 40, 59, 51, 50, 75)_{\adfsmFJj}$,\adfsplit
$(13, 41, 60, 52, 51, 76)_{\adfsmFJj}$,
$(8, 33, 45, 12, 9, 1)_{\adfsmFJj}$,
$(9, 34, 46, 13, 10, 2)_{\adfsmFJj}$,\adfsplit
$(10, 35, 47, 14, 11, 3)_{\adfsmFJj}$,
$(11, 36, 48, 15, 12, 4)_{\adfsmFJj}$,
$(12, 37, 49, 16, 13, 5)_{\adfsmFJj}$,\adfsplit
$(70, 44, 26, 76, 21, 5)_{\adfsmFJj}$,
$(71, 45, 27, 77, 22, 6)_{\adfsmFJj}$,
$(72, 46, 28, 78, 23, 7)_{\adfsmFJj}$,\adfsplit
$(73, 47, 29, 79, 24, 8)_{\adfsmFJj}$,
$(4, 10, 63, 45, 40, 21)_{\adfsmFJj}$,
$(4, 38, 55, 72, 21, 24)_{\adfsmFJj}$,


\adfLgap 
$(34, 26, 35, 58, 77, 15)_{\adfsmFJm}$,
$(64, 25, 55, 38, 78, 30)_{\adfsmFJm}$,
$(43, 32, 80, 37, 73, 27)_{\adfsmFJm}$,\adfsplit
$(53, 22, 19, 3, 41, 16)_{\adfsmFJm}$,
$(83, 57, 31, 79, 23, 36)_{\adfsmFJm}$,
$(12, 73, 72, 51, 33, 20)_{\adfsmFJm}$,\adfsplit
$(75, 14, 45, 54, 0, 27)_{\adfsmFJm}$,
$(83, 45, 47, 66, 34, 25)_{\adfsmFJm}$,
$(75, 19, 44, 72, 63, 52)_{\adfsmFJm}$,\adfsplit
$(37, 23, 84, 33, 38, 35)_{\adfsmFJm}$,
$(1, 66, 25, 2, 16, 60)_{\adfsmFJm}$,
$(13, 11, 33, 34, 56, 30)_{\adfsmFJm}$,\adfsplit
$(27, 67, 76, 79, 10, 64)_{\adfsmFJm}$,
$(25, 41, 81, 69, 30, 53)_{\adfsmFJm}$,
$(1, 39, 26, 14, 8, 40)_{\adfsmFJm}$,\adfsplit
$(56, 46, 42, 19, 72, 1)_{\adfsmFJm}$,
$(41, 39, 32, 38, 45, 47)_{\adfsmFJm}$,
$(63, 0, 30, 7, 36, 45)_{\adfsmFJm}$,\adfsplit
$(29, 9, 23, 84, 18, 79)_{\adfsmFJm}$,
$(51, 62, 34, 44, 50, 27)_{\adfsmFJm}$,
$(6, 39, 59, 37, 62, 2)_{\adfsmFJm}$,

\adfLgap \noindent by the mapping:
$x \mapsto x + 5 j \adfmod{85}$,
$0 \le j < 17$.

\ADFvfyParStart{(85, ((21, 17, ((85, 5)))), -1)} 

\adfDgap
\noindent{\boldmath $ K_{85}$}~
With the point set $Z_{85}$
 the designs are generated from

\adfLgap 
$(84, 13, 4, 67, 62, 16)_{\adfsmFJf}$,
$(38, 25, 21, 62, 11, 32)_{\adfsmFJf}$,
$(1, 27, 48, 65, 40, 29)_{\adfsmFJf}$,\adfsplit
$(29, 62, 8, 17, 6, 65)_{\adfsmFJf}$,
$(3, 1, 7, 78, 51, 0)_{\adfsmFJf}$,
$(73, 52, 61, 21, 65, 11)_{\adfsmFJf}$,\adfsplit
$(65, 54, 73, 83, 47, 79)_{\adfsmFJf}$,
$(57, 72, 76, 15, 28, 42)_{\adfsmFJf}$,
$(15, 54, 24, 46, 19, 39)_{\adfsmFJf}$,\adfsplit
$(24, 43, 7, 63, 59, 26)_{\adfsmFJf}$,
$(32, 31, 48, 26, 23, 64)_{\adfsmFJf}$,
$(38, 39, 33, 28, 76, 27)_{\adfsmFJf}$,\adfsplit
$(79, 42, 52, 82, 24, 17)_{\adfsmFJf}$,
$(37, 42, 24, 10, 38, 6)_{\adfsmFJf}$,
$(19, 82, 5, 70, 65, 73)_{\adfsmFJf}$,\adfsplit
$(71, 30, 65, 40, 4, 10)_{\adfsmFJf}$,
$(64, 60, 34, 13, 0, 18)_{\adfsmFJf}$,


\adfLgap 
$(74, 20, 6, 57, 17, 75)_{\adfsmFJh}$,
$(61, 72, 15, 2, 46, 81)_{\adfsmFJh}$,
$(56, 79, 68, 77, 28, 7)_{\adfsmFJh}$,\adfsplit
$(66, 36, 23, 34, 30, 44)_{\adfsmFJh}$,
$(47, 74, 79, 63, 14, 44)_{\adfsmFJh}$,
$(50, 28, 59, 48, 11, 54)_{\adfsmFJh}$,\adfsplit
$(84, 70, 8, 21, 26, 47)_{\adfsmFJh}$,
$(2, 34, 51, 76, 59, 77)_{\adfsmFJh}$,
$(72, 59, 20, 17, 36, 79)_{\adfsmFJh}$,\adfsplit
$(63, 60, 73, 0, 59, 17)_{\adfsmFJh}$,
$(22, 48, 25, 52, 73, 43)_{\adfsmFJh}$,
$(43, 81, 3, 19, 28, 9)_{\adfsmFJh}$,\adfsplit
$(30, 4, 61, 50, 77, 37)_{\adfsmFJh}$,
$(46, 45, 80, 17, 1, 69)_{\adfsmFJh}$,
$(25, 78, 7, 6, 44, 75)_{\adfsmFJh}$,\adfsplit
$(23, 57, 51, 64, 53, 49)_{\adfsmFJh}$,
$(21, 14, 22, 42, 51, 60)_{\adfsmFJh}$,

\adfLgap \noindent by the mapping:
$x \mapsto x + 4 j \adfmod{84}$ for $x < 84$,
$84 \mapsto 84$,
$0 \le j < 21$.

\ADFvfyParStart{(85, ((17, 21, ((84, 4), (1, 1)))), -1)} 

\adfDgap
\noindent{\boldmath $ K_{96}$}~
With the point set $Z_{96}$
 the designs are generated from

\adfLgap 
$(4, 54, 77, 84, 3, 13)_{\adfsmFJc}$,
$(29, 39, 15, 23, 20, 95)_{\adfsmFJc}$,
$(21, 32, 63, 36, 1, 16)_{\adfsmFJc}$,\adfsplit
$(21, 19, 27, 33, 77, 58)_{\adfsmFJc}$,
$(65, 37, 52, 48, 38, 39)_{\adfsmFJc}$,
$(0, 81, 21, 34, 62, 72)_{\adfsmFJc}$,\adfsplit
$(73, 27, 51, 37, 15, 68)_{\adfsmFJc}$,
$(60, 71, 2, 61, 20, 18)_{\adfsmFJc}$,
$(30, 71, 8, 63, 73, 10)_{\adfsmFJc}$,\adfsplit
$(11, 29, 34, 83, 8, 4)_{\adfsmFJc}$,
$(54, 66, 75, 18, 34, 5)_{\adfsmFJc}$,
$(91, 45, 63, 60, 39, 36)_{\adfsmFJc}$,\adfsplit
$(45, 43, 67, 72, 15, 20)_{\adfsmFJc}$,
$(40, 63, 95, 72, 56, 42)_{\adfsmFJc}$,
$(9, 76, 93, 5, 10, 60)_{\adfsmFJc}$,\adfsplit
$(79, 48, 82, 3, 44, 17)_{\adfsmFJc}$,
$(73, 16, 17, 93, 10, 13)_{\adfsmFJc}$,
$(90, 43, 86, 56, 30, 80)_{\adfsmFJc}$,\adfsplit
$(42, 51, 89, 25, 7, 90)_{\adfsmFJc}$,
$(63, 73, 22, 48, 34, 29)_{\adfsmFJc}$,
$(65, 9, 20, 49, 62, 46)_{\adfsmFJc}$,\adfsplit
$(36, 39, 82, 61, 7, 9)_{\adfsmFJc}$,
$(92, 94, 27, 35, 67, 74)_{\adfsmFJc}$,
$(4, 21, 66, 74, 82, 41)_{\adfsmFJc}$,


\adfLgap 
$(68, 29, 46, 53, 88, 75)_{\adfsmFJj}$,
$(56, 7, 50, 82, 25, 55)_{\adfsmFJj}$,
$(80, 23, 32, 43, 85, 8)_{\adfsmFJj}$,\adfsplit
$(47, 35, 40, 71, 38, 30)_{\adfsmFJj}$,
$(57, 6, 74, 95, 15, 16)_{\adfsmFJj}$,
$(43, 2, 49, 93, 30, 48)_{\adfsmFJj}$,\adfsplit
$(74, 2, 28, 31, 83, 40)_{\adfsmFJj}$,
$(83, 22, 35, 81, 4, 95)_{\adfsmFJj}$,
$(19, 61, 77, 82, 42, 7)_{\adfsmFJj}$,\adfsplit
$(85, 8, 29, 82, 20, 11)_{\adfsmFJj}$,
$(17, 26, 56, 73, 0, 62)_{\adfsmFJj}$,
$(36, 4, 69, 70, 12, 51)_{\adfsmFJj}$,\adfsplit
$(81, 74, 77, 79, 5, 86)_{\adfsmFJj}$,
$(14, 33, 71, 85, 74, 20)_{\adfsmFJj}$,
$(44, 22, 28, 87, 0, 57)_{\adfsmFJj}$,\adfsplit
$(54, 0, 4, 68, 55, 40)_{\adfsmFJj}$,
$(14, 74, 83, 94, 67, 54)_{\adfsmFJj}$,
$(13, 62, 76, 77, 87, 21)_{\adfsmFJj}$,\adfsplit
$(15, 5, 65, 81, 85, 17)_{\adfsmFJj}$,
$(78, 11, 51, 88, 17, 6)_{\adfsmFJj}$,
$(69, 7, 8, 73, 74, 94)_{\adfsmFJj}$,\adfsplit
$(21, 34, 44, 56, 81, 66)_{\adfsmFJj}$,
$(86, 3, 22, 73, 90, 75)_{\adfsmFJj}$,
$(36, 56, 64, 80, 15, 26)_{\adfsmFJj}$,


\adfLgap 
$(95, 79, 43, 32, 10, 31)_{\adfsmFJm}$,
$(92, 23, 36, 64, 1, 10)_{\adfsmFJm}$,
$(76, 26, 66, 85, 83, 81)_{\adfsmFJm}$,\adfsplit
$(43, 73, 67, 62, 65, 46)_{\adfsmFJm}$,
$(46, 3, 15, 85, 70, 60)_{\adfsmFJm}$,
$(83, 17, 78, 36, 85, 6)_{\adfsmFJm}$,\adfsplit
$(16, 52, 41, 82, 14, 69)_{\adfsmFJm}$,
$(4, 61, 55, 41, 93, 20)_{\adfsmFJm}$,
$(14, 29, 36, 91, 73, 78)_{\adfsmFJm}$,\adfsplit
$(68, 52, 66, 53, 5, 22)_{\adfsmFJm}$,
$(55, 40, 60, 52, 13, 61)_{\adfsmFJm}$,
$(10, 28, 0, 94, 67, 42)_{\adfsmFJm}$,\adfsplit
$(28, 19, 63, 54, 40, 53)_{\adfsmFJm}$,
$(33, 41, 94, 89, 70, 78)_{\adfsmFJm}$,
$(32, 15, 22, 77, 14, 54)_{\adfsmFJm}$,\adfsplit
$(29, 24, 2, 16, 20, 60)_{\adfsmFJm}$,
$(13, 57, 17, 59, 22, 34)_{\adfsmFJm}$,
$(22, 34, 91, 92, 58, 37)_{\adfsmFJm}$,\adfsplit
$(22, 16, 19, 2, 65, 64)_{\adfsmFJm}$,
$(0, 76, 72, 87, 65, 41)_{\adfsmFJm}$,
$(18, 12, 73, 93, 46, 52)_{\adfsmFJm}$,\adfsplit
$(79, 1, 93, 72, 89, 50)_{\adfsmFJm}$,
$(65, 48, 59, 45, 29, 79)_{\adfsmFJm}$,
$(16, 59, 0, 71, 19, 39)_{\adfsmFJm}$,

\adfLgap \noindent by the mapping:
$x \mapsto x + 5 j \adfmod{95}$ for $x < 95$,
$95 \mapsto 95$,
$0 \le j < 19$.

\ADFvfyParStart{(96, ((24, 19, ((95, 5), (1, 1)))), -1)} 

\adfDgap
\noindent{\boldmath $ K_{105}$}~
With the point set $Z_{105}$
 the designs are generated from

\adfLgap 
$(43, 50, 99, 41, 7, 16)_{\adfsmFJh}$,
$(71, 83, 90, 76, 12, 30)_{\adfsmFJh}$,
$(33, 17, 35, 91, 9, 24)_{\adfsmFJh}$,\adfsplit
$(19, 102, 74, 91, 97, 55)_{\adfsmFJh}$,
$(38, 69, 41, 60, 49, 27)_{\adfsmFJh}$,
$(78, 44, 96, 50, 102, 1)_{\adfsmFJh}$,\adfsplit
$(26, 15, 41, 52, 58, 3)_{\adfsmFJh}$,
$(24, 17, 26, 25, 92, 84)_{\adfsmFJh}$,
$(72, 12, 82, 65, 31, 95)_{\adfsmFJh}$,\adfsplit
$(10, 19, 15, 35, 50, 7)_{\adfsmFJh}$,
$(43, 63, 60, 8, 76, 59)_{\adfsmFJh}$,
$(21, 15, 39, 4, 66, 88)_{\adfsmFJh}$,\adfsplit
$(66, 103, 88, 3, 62, 17)_{\adfsmFJh}$,
$(60, 101, 66, 94, 36, 91)_{\adfsmFJh}$,
$(21, 18, 78, 84, 10, 71)_{\adfsmFJh}$,\adfsplit
$(4, 32, 92, 63, 16, 23)_{\adfsmFJh}$,
$(42, 78, 82, 103, 57, 41)_{\adfsmFJh}$,
$(99, 93, 15, 55, 27, 8)_{\adfsmFJh}$,\adfsplit
$(3, 4, 45, 74, 95, 44)_{\adfsmFJh}$,
$(98, 84, 1, 23, 69, 22)_{\adfsmFJh}$,
$(104, 36, 72, 50, 94, 65)_{\adfsmFJh}$,\adfsplit
$(7, 89, 92, 58, 94, 64)_{\adfsmFJh}$,
$(27, 37, 51, 43, 71, 100)_{\adfsmFJh}$,
$(70, 5, 72, 71, 22, 50)_{\adfsmFJh}$,\adfsplit
$(26, 48, 38, 30, 74, 49)_{\adfsmFJh}$,
$(50, 27, 79, 73, 6, 40)_{\adfsmFJh}$,


\adfLgap 
$(45, 104, 40, 55, 4, 48)_{\adfsmFJm}$,
$(26, 31, 104, 1, 29, 63)_{\adfsmFJm}$,
$(3, 19, 47, 80, 54, 59)_{\adfsmFJm}$,\adfsplit
$(98, 27, 5, 41, 11, 96)_{\adfsmFJm}$,
$(17, 66, 99, 57, 38, 49)_{\adfsmFJm}$,
$(19, 33, 93, 58, 57, 53)_{\adfsmFJm}$,\adfsplit
$(47, 18, 42, 63, 72, 46)_{\adfsmFJm}$,
$(10, 91, 71, 80, 75, 31)_{\adfsmFJm}$,
$(13, 51, 49, 96, 77, 65)_{\adfsmFJm}$,\adfsplit
$(22, 64, 51, 95, 70, 57)_{\adfsmFJm}$,
$(70, 3, 11, 52, 82, 19)_{\adfsmFJm}$,
$(104, 21, 102, 79, 75, 56)_{\adfsmFJm}$,\adfsplit
$(13, 98, 34, 71, 12, 52)_{\adfsmFJm}$,
$(90, 56, 15, 30, 63, 7)_{\adfsmFJm}$,
$(22, 39, 69, 79, 76, 42)_{\adfsmFJm}$,\adfsplit
$(33, 79, 10, 96, 34, 83)_{\adfsmFJm}$,
$(54, 20, 1, 45, 11, 64)_{\adfsmFJm}$,
$(87, 6, 70, 76, 42, 44)_{\adfsmFJm}$,\adfsplit
$(87, 34, 26, 21, 46, 90)_{\adfsmFJm}$,
$(72, 83, 68, 51, 25, 82)_{\adfsmFJm}$,
$(53, 20, 59, 86, 62, 83)_{\adfsmFJm}$,\adfsplit
$(58, 60, 81, 27, 80, 73)_{\adfsmFJm}$,
$(30, 92, 69, 85, 98, 43)_{\adfsmFJm}$,
$(100, 13, 99, 83, 9, 71)_{\adfsmFJm}$,\adfsplit
$(59, 41, 68, 26, 79, 78)_{\adfsmFJm}$,
$(70, 72, 2, 59, 57, 55)_{\adfsmFJm}$,

\adfLgap \noindent by the mapping:
$x \mapsto x + 5 j \adfmod{105}$,
$0 \le j < 21$.

\ADFvfyParStart{(105, ((26, 21, ((105, 5)))), -1)} 

\adfDgap
\noindent{\boldmath $ K_{116}$}~
With the point set $Z_{116}$
 the designs are generated from

\adfLgap 
$(6, 22, 67, 112, 87, 107)_{\adfsmFJc}$,
$(60, 17, 107, 16, 114, 102)_{\adfsmFJc}$,
$(19, 81, 13, 21, 23, 112)_{\adfsmFJc}$,\adfsplit
$(92, 3, 25, 60, 55, 43)_{\adfsmFJc}$,
$(37, 74, 82, 26, 110, 107)_{\adfsmFJc}$,
$(61, 58, 79, 68, 107, 108)_{\adfsmFJc}$,\adfsplit
$(81, 86, 3, 16, 66, 95)_{\adfsmFJc}$,
$(82, 88, 58, 33, 83, 24)_{\adfsmFJc}$,
$(84, 62, 96, 99, 114, 4)_{\adfsmFJc}$,\adfsplit
$(26, 35, 53, 22, 9, 14)_{\adfsmFJc}$,
$(7, 57, 5, 37, 24, 104)_{\adfsmFJc}$,
$(103, 1, 23, 79, 62, 31)_{\adfsmFJc}$,\adfsplit
$(4, 102, 52, 25, 91, 41)_{\adfsmFJc}$,
$(12, 10, 50, 113, 64, 55)_{\adfsmFJc}$,
$(70, 87, 93, 96, 114, 13)_{\adfsmFJc}$,\adfsplit
$(49, 40, 91, 114, 7, 2)_{\adfsmFJc}$,
$(63, 56, 15, 97, 55, 58)_{\adfsmFJc}$,
$(82, 36, 53, 81, 76, 68)_{\adfsmFJc}$,\adfsplit
$(111, 61, 7, 26, 54, 95)_{\adfsmFJc}$,
$(112, 84, 92, 53, 108, 114)_{\adfsmFJc}$,
$(57, 73, 97, 61, 44, 52)_{\adfsmFJc}$,\adfsplit
$(59, 80, 102, 69, 24, 28)_{\adfsmFJc}$,
$(96, 46, 14, 5, 43, 35)_{\adfsmFJc}$,


\adfLgap 
$(36, 85, 107, 100, 48, 52)_{\adfsmFJf}$,
$(77, 7, 89, 93, 69, 62)_{\adfsmFJf}$,
$(38, 108, 70, 42, 58, 37)_{\adfsmFJf}$,\adfsplit
$(58, 103, 49, 61, 112, 113)_{\adfsmFJf}$,
$(43, 85, 56, 7, 10, 110)_{\adfsmFJf}$,
$(64, 30, 3, 28, 103, 107)_{\adfsmFJf}$,\adfsplit
$(85, 76, 47, 37, 25, 44)_{\adfsmFJf}$,
$(34, 13, 68, 27, 112, 0)_{\adfsmFJf}$,
$(9, 96, 7, 97, 106, 94)_{\adfsmFJf}$,\adfsplit
$(20, 42, 75, 55, 15, 6)_{\adfsmFJf}$,
$(93, 11, 42, 102, 34, 63)_{\adfsmFJf}$,
$(16, 20, 4, 9, 10, 101)_{\adfsmFJf}$,\adfsplit
$(43, 2, 107, 42, 83, 41)_{\adfsmFJf}$,
$(88, 59, 32, 62, 12, 50)_{\adfsmFJf}$,
$(23, 19, 3, 37, 40, 82)_{\adfsmFJf}$,\adfsplit
$(23, 12, 96, 36, 31, 95)_{\adfsmFJf}$,
$(97, 34, 57, 17, 102, 50)_{\adfsmFJf}$,
$(0, 8, 99, 71, 31, 94)_{\adfsmFJf}$,\adfsplit
$(88, 115, 58, 109, 45, 14)_{\adfsmFJf}$,
$(101, 69, 65, 43, 63, 8)_{\adfsmFJf}$,
$(18, 85, 88, 60, 61, 103)_{\adfsmFJf}$,\adfsplit
$(66, 37, 98, 84, 14, 40)_{\adfsmFJf}$,
$(65, 30, 43, 54, 11, 109)_{\adfsmFJf}$,


\adfLgap 
$(55, 54, 88, 10, 21, 100)_{\adfsmFJh}$,
$(20, 90, 94, 71, 73, 24)_{\adfsmFJh}$,
$(32, 9, 97, 14, 104, 62)_{\adfsmFJh}$,\adfsplit
$(6, 90, 15, 81, 61, 34)_{\adfsmFJh}$,
$(113, 51, 63, 54, 28, 34)_{\adfsmFJh}$,
$(59, 72, 96, 63, 19, 94)_{\adfsmFJh}$,\adfsplit
$(115, 80, 6, 78, 109, 114)_{\adfsmFJh}$,
$(65, 57, 21, 58, 23, 46)_{\adfsmFJh}$,
$(49, 88, 26, 47, 32, 29)_{\adfsmFJh}$,\adfsplit
$(80, 45, 2, 14, 42, 107)_{\adfsmFJh}$,
$(98, 105, 53, 8, 80, 3)_{\adfsmFJh}$,
$(34, 81, 107, 92, 20, 104)_{\adfsmFJh}$,\adfsplit
$(21, 14, 75, 98, 3, 71)_{\adfsmFJh}$,
$(60, 17, 88, 59, 96, 11)_{\adfsmFJh}$,
$(69, 41, 109, 11, 31, 28)_{\adfsmFJh}$,\adfsplit
$(114, 27, 8, 90, 14, 3)_{\adfsmFJh}$,
$(64, 72, 69, 51, 42, 73)_{\adfsmFJh}$,
$(90, 95, 26, 39, 75, 7)_{\adfsmFJh}$,\adfsplit
$(62, 5, 37, 52, 97, 54)_{\adfsmFJh}$,
$(109, 19, 58, 3, 83, 76)_{\adfsmFJh}$,
$(28, 29, 68, 59, 14, 13)_{\adfsmFJh}$,\adfsplit
$(105, 36, 109, 91, 1, 11)_{\adfsmFJh}$,
$(8, 4, 11, 104, 19, 72)_{\adfsmFJh}$,


\adfLgap 
$(21, 2, 14, 110, 56, 42)_{\adfsmFJj}$,
$(24, 45, 56, 104, 23, 37)_{\adfsmFJj}$,
$(18, 13, 77, 95, 28, 34)_{\adfsmFJj}$,\adfsplit
$(87, 35, 45, 110, 62, 15)_{\adfsmFJj}$,
$(96, 7, 73, 106, 47, 101)_{\adfsmFJj}$,
$(111, 30, 45, 105, 88, 115)_{\adfsmFJj}$,\adfsplit
$(61, 16, 47, 50, 79, 28)_{\adfsmFJj}$,
$(52, 27, 56, 107, 95, 53)_{\adfsmFJj}$,
$(105, 26, 29, 82, 114, 5)_{\adfsmFJj}$,\adfsplit
$(54, 25, 39, 109, 23, 89)_{\adfsmFJj}$,
$(17, 18, 24, 62, 7, 53)_{\adfsmFJj}$,
$(74, 1, 78, 92, 108, 111)_{\adfsmFJj}$,\adfsplit
$(40, 92, 99, 100, 2, 101)_{\adfsmFJj}$,
$(86, 15, 36, 56, 12, 80)_{\adfsmFJj}$,
$(14, 54, 78, 79, 64, 67)_{\adfsmFJj}$,\adfsplit
$(69, 27, 73, 100, 6, 32)_{\adfsmFJj}$,
$(36, 14, 27, 47, 19, 38)_{\adfsmFJj}$,
$(6, 32, 97, 109, 19, 113)_{\adfsmFJj}$,\adfsplit
$(88, 43, 50, 86, 19, 51)_{\adfsmFJj}$,
$(29, 48, 77, 101, 71, 87)_{\adfsmFJj}$,
$(38, 72, 84, 112, 37, 108)_{\adfsmFJj}$,\adfsplit
$(39, 29, 37, 57, 91, 44)_{\adfsmFJj}$,
$(39, 30, 77, 99, 80, 108)_{\adfsmFJj}$,


\adfLgap 
$(87, 74, 82, 30, 88, 84)_{\adfsmFJm}$,
$(115, 106, 19, 43, 30, 21)_{\adfsmFJm}$,
$(75, 17, 46, 11, 19, 115)_{\adfsmFJm}$,\adfsplit
$(81, 76, 104, 8, 14, 85)_{\adfsmFJm}$,
$(56, 106, 86, 32, 72, 97)_{\adfsmFJm}$,
$(91, 112, 44, 67, 107, 0)_{\adfsmFJm}$,\adfsplit
$(39, 57, 82, 32, 67, 14)_{\adfsmFJm}$,
$(49, 4, 108, 46, 74, 30)_{\adfsmFJm}$,
$(101, 100, 68, 78, 37, 41)_{\adfsmFJm}$,\adfsplit
$(15, 93, 66, 63, 98, 49)_{\adfsmFJm}$,
$(6, 7, 37, 115, 69, 109)_{\adfsmFJm}$,
$(112, 48, 28, 108, 19, 71)_{\adfsmFJm}$,\adfsplit
$(4, 80, 59, 107, 10, 24)_{\adfsmFJm}$,
$(18, 105, 113, 79, 0, 17)_{\adfsmFJm}$,
$(83, 106, 52, 46, 42, 44)_{\adfsmFJm}$,\adfsplit
$(66, 47, 73, 21, 30, 101)_{\adfsmFJm}$,
$(114, 28, 37, 43, 49, 22)_{\adfsmFJm}$,
$(112, 95, 77, 45, 84, 8)_{\adfsmFJm}$,\adfsplit
$(46, 42, 24, 7, 34, 104)_{\adfsmFJm}$,
$(44, 103, 95, 29, 61, 105)_{\adfsmFJm}$,
$(101, 85, 18, 72, 15, 46)_{\adfsmFJm}$,\adfsplit
$(67, 113, 81, 63, 105, 99)_{\adfsmFJm}$,
$(37, 94, 10, 28, 85, 106)_{\adfsmFJm}$,

\adfLgap \noindent by the mapping:
$x \mapsto x + 4 j \adfmod{116}$,
$0 \le j < 29$.

\ADFvfyParStart{(116, ((23, 29, ((116, 4)))), -1)} 

\adfDgap
\noindent{\boldmath $ K_{136}$}~
With the point set $Z_{136}$
 the design is generated from

\adfLgap 
$(135, 131, 57, 134, 130, 58)_{\adfsmFJm}$,
$(31, 15, 132, 9, 14, 71)_{\adfsmFJm}$,
$(68, 29, 41, 130, 63, 110)_{\adfsmFJm}$,\adfsplit
$(9, 69, 40, 102, 33, 81)_{\adfsmFJm}$,
$(54, 92, 122, 110, 99, 57)_{\adfsmFJm}$,
$(114, 60, 63, 71, 49, 23)_{\adfsmFJm}$,\adfsplit
$(71, 82, 19, 45, 121, 88)_{\adfsmFJm}$,
$(112, 131, 23, 100, 103, 83)_{\adfsmFJm}$,
$(121, 21, 20, 72, 58, 126)_{\adfsmFJm}$,\adfsplit
$(53, 92, 90, 28, 46, 96)_{\adfsmFJm}$,
$(50, 103, 68, 1, 58, 115)_{\adfsmFJm}$,
$(122, 19, 7, 62, 28, 12)_{\adfsmFJm}$,\adfsplit
$(12, 63, 2, 84, 133, 15)_{\adfsmFJm}$,
$(95, 37, 17, 81, 120, 74)_{\adfsmFJm}$,
$(61, 4, 125, 109, 65, 89)_{\adfsmFJm}$,\adfsplit
$(54, 56, 127, 94, 109, 53)_{\adfsmFJm}$,
$(122, 34, 105, 63, 82, 14)_{\adfsmFJm}$,
$(28, 78, 99, 109, 31, 85)_{\adfsmFJm}$,\adfsplit
$(15, 30, 44, 79, 70, 54)_{\adfsmFJm}$,
$(78, 46, 86, 104, 105, 96)_{\adfsmFJm}$,
$(78, 42, 134, 12, 97, 4)_{\adfsmFJm}$,\adfsplit
$(55, 14, 68, 97, 65, 91)_{\adfsmFJm}$,
$(17, 41, 86, 32, 98, 16)_{\adfsmFJm}$,
$(70, 97, 54, 62, 48, 31)_{\adfsmFJm}$,\adfsplit
$(18, 3, 31, 14, 56, 2)_{\adfsmFJm}$,
$(45, 10, 103, 40, 73, 47)_{\adfsmFJm}$,
$(17, 31, 81, 76, 105, 39)_{\adfsmFJm}$,\adfsplit
$(75, 98, 73, 30, 21, 97)_{\adfsmFJm}$,
$(125, 38, 49, 73, 36, 40)_{\adfsmFJm}$,
$(15, 118, 112, 107, 3, 71)_{\adfsmFJm}$,\adfsplit
$(134, 132, 12, 88, 45, 57)_{\adfsmFJm}$,
$(21, 36, 17, 9, 124, 126)_{\adfsmFJm}$,
$(96, 17, 116, 45, 23, 36)_{\adfsmFJm}$,\adfsplit
$(84, 69, 2, 119, 65, 103)_{\adfsmFJm}$,

\adfLgap \noindent by the mapping:
$x \mapsto x + 5 j \adfmod{135}$ for $x < 135$,
$135 \mapsto 135$,
$0 \le j < 27$.

\ADFvfyParStart{(136, ((34, 27, ((135, 5), (1, 1)))), -1)} 

\adfDgap
\noindent{\boldmath $ K_{156}$}~
With the point set $Z_{156}$
 the designs are generated from

\adfLgap 
$(53, 2, 27, 106, 144, 41)_{\adfsmFJc}$,
$(141, 49, 78, 107, 82, 72)_{\adfsmFJc}$,
$(102, 151, 20, 71, 103, 134)_{\adfsmFJc}$,\adfsplit
$(64, 76, 111, 140, 59, 155)_{\adfsmFJc}$,
$(108, 4, 22, 29, 42, 149)_{\adfsmFJc}$,
$(67, 46, 111, 138, 10, 144)_{\adfsmFJc}$,\adfsplit
$(21, 43, 94, 103, 131, 81)_{\adfsmFJc}$,
$(72, 120, 128, 81, 150, 52)_{\adfsmFJc}$,
$(10, 91, 16, 150, 79, 112)_{\adfsmFJc}$,\adfsplit
$(126, 141, 63, 113, 87, 137)_{\adfsmFJc}$,
$(140, 29, 44, 141, 55, 145)_{\adfsmFJc}$,
$(5, 134, 147, 6, 94, 104)_{\adfsmFJc}$,\adfsplit
$(11, 101, 137, 151, 13, 31)_{\adfsmFJc}$,
$(129, 134, 136, 8, 132, 137)_{\adfsmFJc}$,
$(102, 68, 92, 99, 35, 54)_{\adfsmFJc}$,\adfsplit
$(105, 10, 33, 3, 57, 115)_{\adfsmFJc}$,
$(132, 122, 50, 111, 64, 29)_{\adfsmFJc}$,
$(99, 94, 118, 37, 144, 56)_{\adfsmFJc}$,\adfsplit
$(9, 85, 126, 30, 134, 76)_{\adfsmFJc}$,
$(72, 144, 21, 104, 38, 118)_{\adfsmFJc}$,
$(141, 122, 22, 60, 80, 85)_{\adfsmFJc}$,\adfsplit
$(16, 32, 35, 118, 119, 110)_{\adfsmFJc}$,
$(5, 47, 70, 2, 103, 16)_{\adfsmFJc}$,
$(29, 19, 47, 148, 60, 98)_{\adfsmFJc}$,\adfsplit
$(149, 131, 153, 132, 67, 106)_{\adfsmFJc}$,
$(32, 45, 91, 152, 5, 23)_{\adfsmFJc}$,
$(66, 0, 101, 110, 107, 54)_{\adfsmFJc}$,\adfsplit
$(28, 10, 121, 141, 140, 43)_{\adfsmFJc}$,
$(146, 140, 69, 7, 103, 66)_{\adfsmFJc}$,
$(141, 103, 135, 99, 143, 1)_{\adfsmFJc}$,\adfsplit
$(19, 108, 137, 66, 85, 54)_{\adfsmFJc}$,


\adfLgap 
$(81, 5, 103, 136, 129, 140)_{\adfsmFJj}$,
$(104, 40, 52, 83, 43, 4)_{\adfsmFJj}$,
$(27, 7, 86, 99, 111, 52)_{\adfsmFJj}$,\adfsplit
$(48, 103, 121, 144, 38, 97)_{\adfsmFJj}$,
$(112, 9, 73, 152, 130, 139)_{\adfsmFJj}$,
$(151, 31, 109, 127, 150, 117)_{\adfsmFJj}$,\adfsplit
$(105, 120, 129, 142, 143, 116)_{\adfsmFJj}$,
$(100, 20, 141, 142, 131, 147)_{\adfsmFJj}$,
$(99, 11, 45, 61, 16, 13)_{\adfsmFJj}$,\adfsplit
$(51, 125, 127, 128, 70, 120)_{\adfsmFJj}$,
$(6, 32, 49, 142, 140, 38)_{\adfsmFJj}$,
$(112, 0, 63, 102, 28, 151)_{\adfsmFJj}$,\adfsplit
$(96, 12, 94, 101, 45, 103)_{\adfsmFJj}$,
$(80, 112, 146, 149, 61, 101)_{\adfsmFJj}$,
$(104, 16, 54, 73, 35, 41)_{\adfsmFJj}$,\adfsplit
$(94, 16, 40, 115, 32, 0)_{\adfsmFJj}$,
$(3, 16, 97, 155, 2, 32)_{\adfsmFJj}$,
$(144, 106, 150, 154, 121, 148)_{\adfsmFJj}$,\adfsplit
$(59, 13, 18, 65, 11, 73)_{\adfsmFJj}$,
$(12, 27, 102, 138, 74, 32)_{\adfsmFJj}$,
$(28, 86, 113, 119, 0, 79)_{\adfsmFJj}$,\adfsplit
$(14, 37, 114, 145, 84, 32)_{\adfsmFJj}$,
$(150, 11, 46, 62, 0, 63)_{\adfsmFJj}$,
$(133, 5, 74, 141, 4, 137)_{\adfsmFJj}$,\adfsplit
$(43, 35, 96, 132, 90, 128)_{\adfsmFJj}$,
$(126, 83, 99, 114, 51, 46)_{\adfsmFJj}$,
$(20, 38, 46, 149, 29, 34)_{\adfsmFJj}$,\adfsplit
$(7, 53, 97, 137, 147, 35)_{\adfsmFJj}$,
$(62, 102, 111, 153, 141, 73)_{\adfsmFJj}$,
$(3, 47, 69, 129, 150, 89)_{\adfsmFJj}$,\adfsplit
$(6, 66, 90, 95, 151, 79)_{\adfsmFJj}$,


\adfLgap 
$(60, 3, 18, 10, 2, 152)_{\adfsmFJm}$,
$(6, 108, 52, 1, 138, 107)_{\adfsmFJm}$,
$(63, 133, 94, 131, 57, 106)_{\adfsmFJm}$,\adfsplit
$(7, 145, 39, 82, 24, 15)_{\adfsmFJm}$,
$(11, 129, 126, 114, 75, 82)_{\adfsmFJm}$,
$(146, 111, 66, 125, 20, 153)_{\adfsmFJm}$,\adfsplit
$(74, 88, 65, 59, 112, 0)_{\adfsmFJm}$,
$(106, 155, 58, 98, 61, 2)_{\adfsmFJm}$,
$(12, 88, 10, 155, 94, 93)_{\adfsmFJm}$,\adfsplit
$(112, 69, 9, 35, 92, 19)_{\adfsmFJm}$,
$(109, 30, 110, 89, 92, 71)_{\adfsmFJm}$,
$(0, 87, 83, 42, 128, 39)_{\adfsmFJm}$,\adfsplit
$(27, 39, 74, 147, 85, 113)_{\adfsmFJm}$,
$(88, 57, 121, 51, 90, 29)_{\adfsmFJm}$,
$(48, 93, 38, 106, 104, 89)_{\adfsmFJm}$,\adfsplit
$(124, 43, 15, 40, 125, 8)_{\adfsmFJm}$,
$(121, 7, 136, 105, 32, 73)_{\adfsmFJm}$,
$(100, 56, 40, 94, 91, 67)_{\adfsmFJm}$,\adfsplit
$(133, 66, 123, 89, 7, 62)_{\adfsmFJm}$,
$(127, 58, 24, 122, 47, 131)_{\adfsmFJm}$,
$(44, 65, 105, 114, 36, 80)_{\adfsmFJm}$,\adfsplit
$(31, 36, 59, 2, 144, 133)_{\adfsmFJm}$,
$(21, 135, 71, 2, 40, 92)_{\adfsmFJm}$,
$(58, 98, 152, 48, 126, 136)_{\adfsmFJm}$,\adfsplit
$(53, 132, 143, 131, 77, 110)_{\adfsmFJm}$,
$(70, 87, 50, 147, 143, 145)_{\adfsmFJm}$,
$(27, 16, 79, 118, 53, 99)_{\adfsmFJm}$,\adfsplit
$(100, 119, 49, 82, 105, 78)_{\adfsmFJm}$,
$(150, 109, 7, 97, 37, 29)_{\adfsmFJm}$,
$(33, 36, 149, 24, 29, 20)_{\adfsmFJm}$,\adfsplit
$(85, 25, 150, 77, 124, 30)_{\adfsmFJm}$,

\adfLgap \noindent by the mapping:
$x \mapsto x + 4 j \adfmod{156}$,
$0 \le j < 39$.

\ADFvfyParStart{(156, ((31, 39, ((156, 4)))), -1)} 

\section{Multipartite graphs}
\label{sec:App-M}
\adfDgap
\noindent{\boldmath $ K_{10,10,10}$}~
With the point set $Z_{30}$
partitioned into
 residue classes modulo $3$,
 the designs are generated from

\adfLgap 
$(0, 1, 3, 5, 14, 23)_{\adfsmFJf}$,


\adfLgap 
$(0, 15, 1, 2, 11, 7)_{\adfsmFJh}$,

\adfLgap \noindent by the mapping:
$x \mapsto x +  j \adfmod{30}$,
$0 \le j < 30$.

\ADFvfyParStart{(30, ((1, 30, ((30, 1)))), ((10, 3)))} 

\adfDgap
\noindent{\boldmath $ K_{10,10,10,10}$}~
With the point set $Z_{40}$
partitioned into
 residue classes modulo $4$,
 the designs are generated from

\adfLgap 
$(0, 9, 22, 35, 19, 21)_{\adfsmFJc}$,
$(13, 10, 11, 36, 20, 19)_{\adfsmFJc}$,
$(36, 1, 19, 30, 34, 35)_{\adfsmFJc}$,


\adfLgap 
$(0, 10, 35, 37, 17, 29)_{\adfsmFJf}$,
$(15, 2, 18, 1, 4, 24)_{\adfsmFJf}$,
$(1, 11, 17, 8, 16, 36)_{\adfsmFJf}$,


\adfLgap 
$(0, 11, 10, 22, 37, 21)_{\adfsmFJh}$,
$(11, 35, 0, 17, 2, 32)_{\adfsmFJh}$,
$(0, 8, 13, 14, 15, 31)_{\adfsmFJh}$,


\adfLgap 
$(0, 2, 7, 33, 24, 37)_{\adfsmFJj}$,
$(13, 3, 14, 28, 9, 31)_{\adfsmFJj}$,
$(0, 1, 3, 30, 24, 27)_{\adfsmFJj}$,


\adfLgap 
$(0, 31, 13, 14, 2, 7)_{\adfsmFJm}$,
$(18, 27, 3, 8, 17, 5)_{\adfsmFJm}$,
$(0, 1, 18, 6, 23, 21)_{\adfsmFJm}$,

\adfLgap \noindent by the mapping:
$x \mapsto x + 2 j \adfmod{40}$,
$0 \le j < 20$.

\ADFvfyParStart{(40, ((3, 20, ((40, 2)))), ((10, 4)))} 

\adfDgap
\noindent{\boldmath $ K_{15,15,15,15}$}~
With the point set $Z_{60}$
partitioned into
 residue classes modulo $3$ for $\{0, 1, \dots, 44\}$, and
 $\{45, 46, \dots, 59\}$,
 the designs are generated from

\adfLgap 
$(0, 37, 17, 52, 51, 29)_{\adfsmFJc}$,
$(54, 26, 28, 15, 21, 14)_{\adfsmFJc}$,
$(14, 10, 56, 24, 33, 13)_{\adfsmFJc}$,


\adfLgap 
$(0, 8, 31, 47, 53, 57)_{\adfsmFJf}$,
$(2, 30, 21, 7, 1, 50)_{\adfsmFJf}$,
$(26, 15, 38, 19, 28, 51)_{\adfsmFJf}$,


\adfLgap 
$(0, 31, 46, 14, 35, 49)_{\adfsmFJh}$,
$(42, 39, 5, 26, 56, 49)_{\adfsmFJh}$,
$(11, 29, 4, 10, 6, 53)_{\adfsmFJh}$,


\adfLgap 
$(0, 54, 7, 23, 6, 41)_{\adfsmFJj}$,
$(45, 10, 8, 18, 46, 4)_{\adfsmFJj}$,
$(12, 1, 32, 46, 6, 25)_{\adfsmFJj}$,


\adfLgap 
$(0, 31, 11, 55, 57, 19)_{\adfsmFJm}$,
$(27, 32, 37, 28, 45, 44)_{\adfsmFJm}$,
$(8, 6, 40, 31, 55, 45)_{\adfsmFJm}$,

\adfLgap \noindent by the mapping:
$x \mapsto x +  j \adfmod{45}$ for $x < 45$,
$x \mapsto (x +  j \adfmod{15}) + 45$ for $x \ge 45$,
$0 \le j < 45$.

\ADFvfyParStart{(60, ((3, 45, ((45, 1), (15, 1)))), ((15, 3), (15, 1)))} 

\adfDgap
\noindent{\boldmath $ K_{20,20,20,20}$}~
With the point set $Z_{80}$
partitioned into
 residue classes modulo $4$,
 the designs are generated from

\adfLgap 
$(0, 21, 63, 66, 26, 58)_{\adfsmFJc}$,
$(54, 69, 56, 63, 7, 29)_{\adfsmFJc}$,
$(5, 16, 35, 6, 62, 44)_{\adfsmFJc}$,


\adfLgap 
$(0, 3, 52, 54, 46, 14)_{\adfsmFJf}$,
$(60, 65, 43, 10, 42, 50)_{\adfsmFJf}$,
$(69, 11, 3, 20, 24, 30)_{\adfsmFJf}$,


\adfLgap 
$(0, 48, 18, 73, 69, 54)_{\adfsmFJh}$,
$(26, 77, 67, 39, 40, 72)_{\adfsmFJh}$,
$(24, 10, 1, 7, 59, 9)_{\adfsmFJh}$,


\adfLgap 
$(0, 49, 51, 26, 16, 71)_{\adfsmFJj}$,
$(44, 78, 25, 3, 40, 57)_{\adfsmFJj}$,
$(61, 2, 64, 75, 1, 11)_{\adfsmFJj}$,


\adfLgap 
$(0, 34, 3, 43, 53, 29)_{\adfsmFJm}$,
$(2, 33, 23, 0, 8, 13)_{\adfsmFJm}$,
$(26, 3, 64, 8, 25, 71)_{\adfsmFJm}$,

\adfLgap \noindent by the mapping:
$x \mapsto x +  j \adfmod{80}$,
$0 \le j < 80$.

\ADFvfyParStart{(80, ((3, 80, ((80, 1)))), ((20, 4)))} 

\adfDgap
\noindent{\boldmath $ K_{25,25,25,25}$}~
With the point set $Z_{100}$
partitioned into
 residue classes modulo $3$ for $\{0, 1, \dots, 74\}$, and
 $\{75, 76, \dots, 99\}$,
 the designs are generated from

\adfLgap 
$(0, 56, 52, 83, 94, 25)_{\adfsmFJc}$,
$(94, 10, 48, 26, 20, 3)_{\adfsmFJc}$,
$(33, 93, 20, 31, 40, 80)_{\adfsmFJc}$,\adfsplit
$(49, 66, 81, 5, 17, 57)_{\adfsmFJc}$,
$(0, 34, 79, 29, 74, 86)_{\adfsmFJc}$,


\adfLgap 
$(0, 32, 56, 94, 34, 4)_{\adfsmFJf}$,
$(29, 4, 8, 9, 21, 69)_{\adfsmFJf}$,
$(0, 31, 7, 82, 90, 86)_{\adfsmFJf}$,\adfsplit
$(46, 53, 35, 72, 95, 99)_{\adfsmFJf}$,
$(0, 81, 77, 11, 29, 59)_{\adfsmFJf}$,


\adfLgap 
$(0, 59, 52, 64, 83, 91)_{\adfsmFJh}$,
$(45, 48, 83, 90, 47, 67)_{\adfsmFJh}$,
$(98, 61, 45, 47, 11, 51)_{\adfsmFJh}$,\adfsplit
$(42, 60, 22, 34, 81, 2)_{\adfsmFJh}$,
$(0, 15, 44, 62, 94, 80)_{\adfsmFJh}$,


\adfLgap 
$(0, 93, 71, 52, 18, 86)_{\adfsmFJj}$,
$(78, 19, 14, 30, 81, 70)_{\adfsmFJj}$,
$(13, 87, 42, 56, 55, 82)_{\adfsmFJj}$,\adfsplit
$(59, 33, 61, 88, 53, 34)_{\adfsmFJj}$,
$(0, 10, 17, 88, 48, 35)_{\adfsmFJj}$,


\adfLgap 
$(0, 4, 80, 98, 47, 37)_{\adfsmFJm}$,
$(64, 23, 3, 54, 80, 53)_{\adfsmFJm}$,
$(78, 65, 53, 39, 13, 58)_{\adfsmFJm}$,\adfsplit
$(13, 72, 6, 94, 84, 59)_{\adfsmFJm}$,
$(75, 1, 8, 74, 21, 16)_{\adfsmFJm}$,

\adfLgap \noindent by the mapping:
$x \mapsto x +  j \adfmod{75}$ for $x < 75$,
$x \mapsto (x +  j \adfmod{25}) + 75$ for $x \ge 75$,
$0 \le j < 75$.

\ADFvfyParStart{(100, ((5, 75, ((75, 1), (25, 1)))), ((25, 3), (25, 1)))} 

\adfDgap
\noindent{\boldmath $ K_{5,5,5,9}$}~
With the point set $Z_{24}$
partitioned into
 residue classes modulo $3$ for $\{0, 1, \dots, 14\}$, and
 $\{15, 16, \dots, 23\}$,
 the designs are generated from

\adfLgap 
$(9, 4, 6, 20, 15, 16)_{\adfsmFJf}$,
$(9, 5, 6, 7, 17, 1)_{\adfsmFJf}$,
$(0, 14, 12, 1, 15, 16)_{\adfsmFJf}$,\adfsplit
$(3, 8, 12, 4, 21, 22)_{\adfsmFJf}$,
$(1, 3, 7, 15, 16, 20)_{\adfsmFJf}$,
$(0, 10, 1, 8, 18, 19)_{\adfsmFJf}$,\adfsplit
$(3, 7, 10, 2, 17, 23)_{\adfsmFJf}$,


\adfLgap 
$(16, 21, 0, 6, 2, 11)_{\adfsmFJh}$,
$(2, 8, 1, 12, 23, 19)_{\adfsmFJh}$,
$(15, 17, 12, 8, 14, 10)_{\adfsmFJh}$,\adfsplit
$(15, 19, 0, 3, 4, 13)_{\adfsmFJh}$,
$(0, 20, 1, 7, 5, 8)_{\adfsmFJh}$,
$(16, 17, 1, 4, 9, 5)_{\adfsmFJh}$,\adfsplit
$(1, 4, 3, 12, 14, 21)_{\adfsmFJh}$,

\adfLgap \noindent by the mapping:
$x \mapsto x + 5 j \adfmod{15}$ for $x < 15$,
$x \mapsto (x - 15 + 3 j \adfmod{9}) + 15$ for $x \ge 15$,
$0 \le j < 3$.

\ADFvfyParStart{(24, ((7, 3, ((15, 5), (9, 3)))), ((5, 3), (9, 1)))} 

\adfDgap
\noindent{\boldmath $ K_{10,10,10,15}$}~
With the point set $Z_{45}$
partitioned into
 residue classes modulo $3$ for $\{0, 1, \dots, 29\}$, and
 $\{30, 31, \dots, 44\}$,
 the designs are generated from

\adfLgap 
$(0, 13, 14, 40, 43, 32)_{\adfsmFJc}$,
$(0, 35, 26, 7, 1, 19)_{\adfsmFJc}$,
$(14, 16, 6, 37, 34, 31)_{\adfsmFJc}$,\adfsplit
$(31, 16, 11, 3, 9, 7)_{\adfsmFJc}$,
$(1, 5, 15, 31, 40, 36)_{\adfsmFJc}$,


\adfLgap 
$(0, 14, 7, 33, 36, 39)_{\adfsmFJf}$,
$(31, 23, 32, 15, 19, 18)_{\adfsmFJf}$,
$(20, 35, 36, 9, 16, 19)_{\adfsmFJf}$,\adfsplit
$(37, 12, 21, 4, 7, 19)_{\adfsmFJf}$,
$(0, 11, 8, 1, 28, 34)_{\adfsmFJf}$,


\adfLgap 
$(0, 37, 28, 20, 11, 19)_{\adfsmFJh}$,
$(16, 26, 34, 30, 12, 21)_{\adfsmFJh}$,
$(8, 2, 9, 25, 34, 35)_{\adfsmFJh}$,\adfsplit
$(37, 42, 1, 4, 3, 26)_{\adfsmFJh}$,
$(1, 13, 17, 23, 39, 40)_{\adfsmFJh}$,


\adfLgap 
$(0, 1, 2, 37, 24, 17)_{\adfsmFJj}$,
$(36, 9, 14, 7, 38, 21)_{\adfsmFJj}$,
$(41, 5, 25, 21, 40, 29)_{\adfsmFJj}$,\adfsplit
$(42, 6, 2, 16, 36, 8)_{\adfsmFJj}$,
$(0, 5, 13, 42, 24, 32)_{\adfsmFJj}$,


\adfLgap 
$(0, 2, 13, 7, 38, 33)_{\adfsmFJm}$,
$(12, 23, 8, 35, 36, 1)_{\adfsmFJm}$,
$(5, 19, 13, 15, 40, 35)_{\adfsmFJm}$,\adfsplit
$(22, 0, 6, 17, 39, 30)_{\adfsmFJm}$,
$(0, 1, 20, 29, 41, 35)_{\adfsmFJm}$,

\adfLgap \noindent by the mapping:
$x \mapsto x + 2 j \adfmod{30}$ for $x < 30$,
$x \mapsto (x +  j \adfmod{15}) + 30$ for $x \ge 30$,
$0 \le j < 15$.

\ADFvfyParStart{(45, ((5, 15, ((30, 2), (15, 1)))), ((10, 3), (15, 1)))} 

\adfDgap
\noindent{\boldmath $ K_{3,3,3,3,3}$}~
With the point set $Z_{15}$
partitioned into
 residue classes modulo $5$,
 the design is generated from

\adfLgap 
$(0, 9, 7, 8, 3, 13)_{\adfsmFJf}$,
$(0, 4, 1, 2, 7, 12)_{\adfsmFJf}$,
$(0, 14, 3, 1, 6, 11)_{\adfsmFJf}$,

\adfLgap \noindent by the mapping:
$x \mapsto x + 5 j \adfmod{15}$,
$0 \le j < 3$.

\ADFvfyParStart{(15, ((3, 3, ((15, 5)))), ((3, 5)))} 

\adfDgap
\noindent{\boldmath $ K_{3,3,3,3,3}$}~
With the point set $Z_{15}$
partitioned into
 residue classes modulo $5$,
 the design is generated from

\adfLgap 
$(0, 11, 4, 12, 3, 13)_{\adfsmFJj}$,
$(9, 11, 5, 2, 8, 7)_{\adfsmFJj}$,
$(9, 13, 10, 12, 1, 8)_{\adfsmFJj}$,\adfsplit
$(9, 3, 0, 6, 7, 1)_{\adfsmFJj}$,
$(10, 6, 2, 4, 13, 8)_{\adfsmFJj}$,
$(8, 4, 7, 1, 5, 0)_{\adfsmFJj}$,\adfsplit
$(3, 1, 2, 14, 0, 5)_{\adfsmFJj}$,
$(5, 6, 12, 14, 8, 13)_{\adfsmFJj}$,
$(10, 7, 11, 14, 13, 3)_{\adfsmFJj}$.


\ADFvfyParStart{(15, ((9, 1, ((15, 15)))), ((3, 5)))} 

\adfDgap
\noindent{\boldmath $ K_{5,5,5,5,5}$}~
With the point set $Z_{25}$
partitioned into
 residue classes modulo $5$,
 the designs are generated from

\adfLgap 
$(0, 1, 3, 7, 12, 8)_{\adfsmFJc}$,


\adfLgap 
$(0, 1, 5, 7, 14, 22)_{\adfsmFJf}$,


\adfLgap 
$(0, 5, 1, 2, 18, 11)_{\adfsmFJh}$,


\adfLgap 
$(0, 1, 3, 7, 15, 9)_{\adfsmFJj}$,


\adfLgap 
$(0, 1, 2, 4, 8, 13)_{\adfsmFJm}$,

\adfLgap \noindent by the mapping:
$x \mapsto x +  j \adfmod{25}$,
$0 \le j < 25$.

\ADFvfyParStart{(25, ((1, 25, ((25, 1)))), ((5, 5)))} 

\adfDgap
\noindent{\boldmath $ K_{6,6,6,6,6}$}~
With the point set $Z_{30}$
partitioned into
 residue classes modulo $4$ for $\{0, 1, \dots, 23\}$, and
 $\{24, 25, \dots, 29\}$,
 the designs are generated from

\adfLgap 
$(0, 27, 21, 22, 3, 13)_{\adfsmFJc}$,
$(5, 26, 12, 18, 3, 22)_{\adfsmFJc}$,
$(11, 27, 16, 6, 1, 10)_{\adfsmFJc}$,


\adfLgap 
$(0, 1, 2, 7, 11, 23)_{\adfsmFJf}$,
$(0, 3, 1, 10, 24, 27)_{\adfsmFJf}$,
$(0, 13, 16, 18, 25, 28)_{\adfsmFJf}$,


\adfLgap 
$(0, 8, 19, 5, 6, 28)_{\adfsmFJh}$,
$(21, 12, 25, 22, 15, 26)_{\adfsmFJh}$,
$(1, 9, 0, 11, 29, 18)_{\adfsmFJh}$,


\adfLgap 
$(0, 1, 11, 24, 2, 17)_{\adfsmFJj}$,
$(16, 22, 19, 24, 17, 6)_{\adfsmFJj}$,
$(19, 26, 1, 12, 10, 14)_{\adfsmFJj}$,


\adfLgap 
$(0, 19, 2, 25, 24, 3)_{\adfsmFJm}$,
$(10, 9, 23, 19, 16, 1)_{\adfsmFJm}$,
$(6, 17, 11, 20, 25, 24)_{\adfsmFJm}$,

\adfLgap \noindent by the mapping:
$x \mapsto x + 2 j \adfmod{24}$ for $x < 24$,
$x \mapsto (x +  j \adfmod{6}) + 24$ for $x \ge 24$,
$0 \le j < 12$.

\ADFvfyParStart{(30, ((3, 12, ((24, 2), (6, 1)))), ((6, 4), (6, 1)))} 

\adfDgap
\noindent{\boldmath $ K_{8,8,8,8,8}$}~
With the point set $Z_{40}$
partitioned into
 residue classes modulo $4$ for $\{0, 1, \dots, 31\}$, and
 $\{32, 33, \dots, 39\}$,
 the design is generated from

\adfLgap 
$(0, 1, 2, 7, 32, 29)_{\adfsmFJm}$,
$(32, 3, 4, 14, 21, 23)_{\adfsmFJm}$,

\adfLgap \noindent by the mapping:
$x \mapsto x +  j \adfmod{32}$ for $x < 32$,
$x \mapsto (x +  j \adfmod{8}) + 32$ for $x \ge 32$,
$0 \le j < 32$.

\ADFvfyParStart{(40, ((2, 32, ((32, 1), (8, 1)))), ((8, 4), (8, 1)))} 

\adfDgap
\noindent{\boldmath $ K_{10,10,10,10,10}$}~
With the point set $Z_{50}$
partitioned into
 residue classes modulo $5$,
 the designs are generated from

\adfLgap 
$(0, 48, 34, 27, 26, 39)_{\adfsmFJc}$,
$(8, 7, 20, 11, 26, 25)_{\adfsmFJc}$,


\adfLgap 
$(0, 32, 24, 1, 8, 21)_{\adfsmFJf}$,
$(0, 13, 5, 7, 17, 41)_{\adfsmFJf}$,


\adfLgap 
$(0, 38, 4, 17, 1, 11)_{\adfsmFJh}$,
$(0, 10, 12, 24, 19, 32)_{\adfsmFJh}$,


\adfLgap 
$(0, 44, 12, 41, 45, 7)_{\adfsmFJj}$,
$(0, 8, 22, 24, 35, 19)_{\adfsmFJj}$,


\adfLgap 
$(0, 48, 27, 39, 1, 6)_{\adfsmFJm}$,
$(0, 12, 32, 28, 19, 36)_{\adfsmFJm}$,

\adfLgap \noindent by the mapping:
$x \mapsto x +  j \adfmod{50}$,
$0 \le j < 50$.

\ADFvfyParStart{(50, ((2, 50, ((50, 1)))), ((10, 5)))} 

\adfDgap
\noindent{\boldmath $ K_{21,21,21,21,21}$}~
With the point set $Z_{105}$
partitioned into
 residue classes modulo $5$,
 the design is generated from

\adfLgap 
$(33, 5, 36, 42, 72, 4)_{\adfsmFJc}$,
$(34, 6, 37, 43, 73, 5)_{\adfsmFJc}$,
$(35, 7, 38, 44, 74, 6)_{\adfsmFJc}$,\adfsplit
$(36, 8, 39, 45, 75, 7)_{\adfsmFJc}$,
$(37, 9, 40, 46, 76, 8)_{\adfsmFJc}$,
$(40, 73, 84, 92, 27, 83)_{\adfsmFJc}$,\adfsplit
$(41, 74, 85, 93, 28, 84)_{\adfsmFJc}$,
$(42, 75, 86, 94, 29, 85)_{\adfsmFJc}$,
$(43, 76, 87, 95, 30, 86)_{\adfsmFJc}$,\adfsplit
$(44, 77, 88, 96, 31, 87)_{\adfsmFJc}$,
$(6, 40, 57, 8, 64, 95)_{\adfsmFJc}$,
$(7, 41, 58, 9, 65, 96)_{\adfsmFJc}$,\adfsplit
$(8, 42, 59, 10, 66, 97)_{\adfsmFJc}$,
$(9, 43, 60, 11, 67, 98)_{\adfsmFJc}$,
$(10, 44, 61, 12, 68, 99)_{\adfsmFJc}$,\adfsplit
$(38, 65, 42, 64, 24, 50)_{\adfsmFJc}$,
$(39, 66, 43, 65, 25, 51)_{\adfsmFJc}$,
$(40, 67, 44, 66, 26, 52)_{\adfsmFJc}$,\adfsplit
$(15, 11, 38, 102, 37, 36)_{\adfsmFJc}$,
$(49, 7, 28, 70, 91, 37)_{\adfsmFJc}$,
$(1, 24, 102, 23, 88, 13)_{\adfsmFJc}$,

\adfLgap \noindent by the mapping:
$x \mapsto x + 5 j \adfmod{105}$,
$0 \le j < 21$.

\ADFvfyParStart{(105, ((21, 21, ((105, 5)))), ((21, 5)))} 

\adfDgap
\noindent{\boldmath $ K_{8,8,8,8,3}$}~
With the point set $Z_{35}$
partitioned into
 residue classes modulo $3$ for $\{0, 1, \dots, 23\}$,
 $\{24, 25, \dots, 31\}$, and
 $\{32, 33, 34\}$,
 the design is generated from

\adfLgap 
$(12, 22, 20, 24, 34, 7)_{\adfsmFJm}$,
$(32, 15, 19, 17, 27, 24)_{\adfsmFJm}$,
$(18, 31, 29, 1, 14, 11)_{\adfsmFJm}$,\adfsplit
$(0, 1, 13, 2, 5, 24)_{\adfsmFJm}$,

\adfLgap \noindent by the mapping:
$x \mapsto x + 2 j \adfmod{24}$ for $x < 24$,
$x \mapsto (x + 2 j \adfmod{8}) + 24$ for $24 \le x < 32$,
$x \mapsto (x - 32 +  j \adfmod{3}) + 32$ for $x \ge 32$,
$0 \le j < 12$.

\ADFvfyParStart{(35, ((4, 12, ((24, 2), (8, 2), (3, 1)))), ((8, 3), (8, 1), (3, 1)))} 

\adfDgap
\noindent{\boldmath $ K_{10,10,10,10,15}$}~
With the point set $Z_{55}$
partitioned into
 residue classes modulo $4$ for $\{0, 1, \dots, 39\}$, and
 $\{40, 41, \dots, 54\}$,
 the designs are generated from

\adfLgap 
$(0, 5, 31, 48, 43, 27)_{\adfsmFJc}$,
$(36, 13, 14, 48, 3, 49)_{\adfsmFJc}$,
$(39, 16, 37, 44, 26, 43)_{\adfsmFJc}$,\adfsplit
$(33, 52, 14, 39, 0, 44)_{\adfsmFJc}$,
$(36, 11, 54, 2, 33, 31)_{\adfsmFJc}$,
$(0, 1, 38, 44, 53, 49)_{\adfsmFJc}$,


\adfLgap 
$(0, 40, 35, 13, 1, 5)_{\adfsmFJf}$,
$(42, 13, 38, 28, 24, 36)_{\adfsmFJf}$,
$(32, 39, 23, 53, 50, 45)_{\adfsmFJf}$,\adfsplit
$(37, 11, 24, 18, 2, 45)_{\adfsmFJf}$,
$(5, 50, 9, 6, 18, 34)_{\adfsmFJf}$,
$(0, 21, 26, 23, 43, 46)_{\adfsmFJf}$,


\adfLgap 
$(0, 11, 42, 48, 13, 2)_{\adfsmFJh}$,
$(46, 48, 39, 6, 28, 33)_{\adfsmFJh}$,
$(22, 31, 41, 32, 13, 25)_{\adfsmFJh}$,\adfsplit
$(37, 2, 36, 27, 49, 20)_{\adfsmFJh}$,
$(39, 5, 24, 10, 41, 43)_{\adfsmFJh}$,
$(1, 53, 4, 15, 18, 36)_{\adfsmFJh}$,


\adfLgap 
$(0, 38, 21, 49, 3, 41)_{\adfsmFJj}$,
$(25, 12, 15, 26, 48, 47)_{\adfsmFJj}$,
$(47, 19, 10, 21, 40, 27)_{\adfsmFJj}$,\adfsplit
$(36, 6, 31, 37, 45, 13)_{\adfsmFJj}$,
$(31, 24, 17, 48, 30, 12)_{\adfsmFJj}$,
$(41, 2, 24, 39, 49, 8)_{\adfsmFJj}$,

\adfLgap \noindent by the mapping:
$x \mapsto x + 2 j \adfmod{40}$ for $x < 40$,
$x \mapsto (x - 40 + 3 j \adfmod{15}) + 40$ for $x \ge 40$,
$0 \le j < 20$.

\ADFvfyParStart{(55, ((6, 20, ((40, 2), (15, 3)))), ((10, 4), (15, 1)))} 

\adfDgap
\noindent{\boldmath $ K_{10,10,10,10,15}$}~
With the point set $Z_{55}$
partitioned into
 $\{0, 1, \dots, 14\}$,
 residue classes modulo $3$ for $\{15, 16, \dots, 44\}$, and
 $\{45, 46, \dots, 54\}$,
 the design is generated from

\adfLgap 
$(0, 46, 39, 15, 26, 31)_{\adfsmFJm}$,
$(2, 47, 31, 42, 21, 49)_{\adfsmFJm}$,
$(3, 35, 24, 48, 31, 22)_{\adfsmFJm}$,\adfsplit
$(2, 53, 38, 32, 24, 27)_{\adfsmFJm}$,
$(11, 22, 36, 32, 50, 28)_{\adfsmFJm}$,
$(12, 51, 35, 27, 18, 54)_{\adfsmFJm}$,\adfsplit
$(1, 15, 37, 34, 44, 39)_{\adfsmFJm}$,
$(0, 19, 21, 54, 51, 22)_{\adfsmFJm}$,

\adfLgap \noindent by the mapping:
$x \mapsto x +  j \adfmod{15}$ for $x < 15$,
$x \mapsto (x - 15 + 2 j \adfmod{30}) + 15$ for $15 \le x < 45$,
$x \mapsto (x - 45 + 2 j \adfmod{10}) + 45$ for $x \ge 45$,
$0 \le j < 15$.

\ADFvfyParStart{(55, ((8, 15, ((15, 1), (30, 2), (10, 2)))), ((15, 1), (10, 3), (10, 1)))} 

\adfDgap
\noindent{\boldmath $ K_{10,10,10,10,20}$}~
With the point set $Z_{60}$
partitioned into
 residue classes modulo $4$ for $\{0, 1, \dots, 39\}$, and
 $\{40, 41, \dots, 59\}$,
 the designs are generated from

\adfLgap 
$(48, 18, 27, 28, 24, 34)_{\adfsmFJc}$,
$(31, 6, 4, 46, 43, 51)_{\adfsmFJc}$,
$(25, 56, 28, 2, 6, 8)_{\adfsmFJc}$,\adfsplit
$(18, 29, 31, 41, 56, 54)_{\adfsmFJc}$,
$(46, 17, 36, 27, 3, 16)_{\adfsmFJc}$,
$(0, 1, 35, 49, 57, 5)_{\adfsmFJc}$,\adfsplit
$(0, 15, 33, 46, 50, 29)_{\adfsmFJc}$,


\adfLgap 
$(41, 31, 21, 2, 4, 24)_{\adfsmFJf}$,
$(52, 4, 36, 38, 19, 14)_{\adfsmFJf}$,
$(56, 10, 24, 9, 19, 15)_{\adfsmFJf}$,\adfsplit
$(55, 3, 8, 1, 29, 21)_{\adfsmFJf}$,
$(15, 9, 33, 46, 40, 44)_{\adfsmFJf}$,
$(16, 30, 8, 56, 42, 50)_{\adfsmFJf}$,\adfsplit
$(0, 11, 38, 1, 43, 56)_{\adfsmFJf}$,


\adfLgap 
$(32, 53, 29, 33, 39, 30)_{\adfsmFJh}$,
$(6, 54, 15, 17, 28, 23)_{\adfsmFJh}$,
$(56, 44, 17, 22, 18, 13)_{\adfsmFJh}$,\adfsplit
$(1, 41, 39, 26, 18, 36)_{\adfsmFJh}$,
$(11, 27, 46, 38, 32, 45)_{\adfsmFJh}$,
$(48, 56, 34, 38, 28, 24)_{\adfsmFJh}$,\adfsplit
$(1, 9, 16, 23, 49, 55)_{\adfsmFJh}$,


\adfLgap 
$(21, 20, 52, 38, 7, 8)_{\adfsmFJj}$,
$(50, 31, 2, 13, 49, 29)_{\adfsmFJj}$,
$(45, 14, 29, 35, 47, 20)_{\adfsmFJj}$,\adfsplit
$(25, 15, 52, 30, 36, 27)_{\adfsmFJj}$,
$(48, 3, 4, 6, 53, 21)_{\adfsmFJj}$,
$(47, 16, 21, 30, 51, 6)_{\adfsmFJj}$,\adfsplit
$(0, 7, 30, 33, 42, 41)_{\adfsmFJj}$,


\adfLgap 
$(13, 44, 8, 19, 2, 54)_{\adfsmFJm}$,
$(33, 34, 58, 46, 11, 12)_{\adfsmFJm}$,
$(5, 36, 46, 51, 19, 18)_{\adfsmFJm}$,\adfsplit
$(50, 7, 9, 4, 12, 19)_{\adfsmFJm}$,
$(46, 16, 0, 29, 2, 10)_{\adfsmFJm}$,
$(28, 35, 29, 6, 48, 50)_{\adfsmFJm}$,\adfsplit
$(1, 8, 16, 39, 43, 59)_{\adfsmFJm}$,

\adfLgap \noindent by the mapping:
$x \mapsto x + 2 j \adfmod{40}$ for $x < 40$,
$x \mapsto (x +  j \adfmod{20}) + 40$ for $x \ge 40$,
$0 \le j < 20$.

\ADFvfyParStart{(60, ((7, 20, ((40, 2), (20, 1)))), ((10, 4), (20, 1)))} 

\adfDgap
\noindent{\boldmath $ K_{21,21,21,21,36}$}~
With the point set $Z_{120}$
partitioned into
 residue classes modulo $4$ for $\{0, 1, \dots, 83\}$, and
 $\{84, 85, \dots, 119\}$,
 the design is generated from

\adfLgap 
$(49, 39, 71, 104, 70, 2)_{\adfsmFJm}$,
$(16, 30, 61, 103, 23, 46)_{\adfsmFJm}$,
$(40, 51, 37, 108, 84, 50)_{\adfsmFJm}$,\adfsplit
$(73, 30, 44, 117, 108, 19)_{\adfsmFJm}$,
$(10, 60, 8, 102, 116, 83)_{\adfsmFJm}$,
$(30, 3, 57, 112, 102, 67)_{\adfsmFJm}$,\adfsplit
$(67, 68, 2, 105, 102, 16)_{\adfsmFJm}$,
$(29, 12, 27, 86, 85, 82)_{\adfsmFJm}$,
$(80, 54, 1, 93, 92, 74)_{\adfsmFJm}$,\adfsplit
$(33, 62, 34, 107, 56, 44)_{\adfsmFJm}$,
$(72, 43, 29, 89, 110, 11)_{\adfsmFJm}$,
$(13, 116, 91, 47, 75, 44)_{\adfsmFJm}$,\adfsplit
$(64, 13, 113, 102, 30, 29)_{\adfsmFJm}$,
$(25, 18, 75, 89, 96, 40)_{\adfsmFJm}$,
$(55, 85, 118, 42, 28, 24)_{\adfsmFJm}$,\adfsplit
$(7, 65, 100, 99, 72, 30)_{\adfsmFJm}$,
$(64, 27, 87, 62, 6, 59)_{\adfsmFJm}$,
$(21, 59, 8, 106, 118, 46)_{\adfsmFJm}$,\adfsplit
$(43, 69, 119, 91, 4, 82)_{\adfsmFJm}$,
$(7, 2, 93, 86, 21, 73)_{\adfsmFJm}$,
$(49, 106, 12, 14, 54, 55)_{\adfsmFJm}$,\adfsplit
$(29, 43, 74, 107, 28, 23)_{\adfsmFJm}$,
$(40, 88, 99, 65, 13, 62)_{\adfsmFJm}$,
$(40, 10, 111, 107, 27, 61)_{\adfsmFJm}$,\adfsplit
$(22, 45, 113, 85, 47, 13)_{\adfsmFJm}$,
$(80, 62, 59, 112, 85, 5)_{\adfsmFJm}$,
$(7, 82, 115, 98, 49, 24)_{\adfsmFJm}$,

\adfLgap \noindent by the mapping:
$x \mapsto x + 4 j \adfmod{84}$ for $x < 84$,
$x \mapsto (x - 84 + 12 j \adfmod{36}) + 84$ for $x \ge 84$,
$0 \le j < 21$.

\ADFvfyParStart{(120, ((27, 21, ((84, 4), (36, 12)))), ((21, 4), (36, 1)))} 

\adfDgap
\noindent{\boldmath $ K_{4^{6}}$}~
With the point set $Z_{24}$
partitioned into
 residue classes modulo $6$,
 the designs are generated from

\adfLgap 
$(0, 1, 3, 11, 20, 9)_{\adfsmFJc}$,


\adfLgap 
$(0, 1, 6, 8, 11, 21)_{\adfsmFJf}$,


\adfLgap 
$(0, 1, 2, 5, 11, 8)_{\adfsmFJh}$,


\adfLgap 
$(0, 1, 3, 8, 12, 10)_{\adfsmFJj}$,


\adfLgap 
$(0, 1, 2, 17, 11, 21)_{\adfsmFJm}$,

\adfLgap \noindent by the mapping:
$x \mapsto x +  j \adfmod{24}$,
$0 \le j < 24$.

\ADFvfyParStart{(24, ((1, 24, ((24, 1)))), ((4, 6)))} 

\adfDgap
\noindent{\boldmath $ K_{4,4,4,4,4,7}$}~
With the point set $Z_{27}$
partitioned into
 residue classes modulo $5$ for $\{0, 1, \dots, 19\}$, and
 $\{20, 21, \dots, 26\}$,
 the designs are generated from

\adfLgap 
$(20, 0, 8, 2, 17, 14)_{\adfsmFJc}$,
$(1, 13, 23, 19, 2, 5)_{\adfsmFJc}$,
$(22, 11, 0, 4, 19, 13)_{\adfsmFJc}$,\adfsplit
$(16, 17, 25, 19, 14, 23)_{\adfsmFJc}$,
$(1, 14, 26, 8, 15, 22)_{\adfsmFJc}$,
$(2, 11, 24, 14, 15, 6)_{\adfsmFJc}$,


\adfLgap 
$(25, 13, 0, 12, 19, 2)_{\adfsmFJf}$,
$(10, 24, 7, 14, 4, 13)_{\adfsmFJf}$,
$(8, 5, 9, 1, 21, 24)_{\adfsmFJf}$,\adfsplit
$(2, 11, 19, 20, 3, 24)_{\adfsmFJf}$,
$(23, 8, 1, 2, 12, 19)_{\adfsmFJf}$,
$(3, 26, 14, 1, 6, 16)_{\adfsmFJf}$,


\adfLgap 
$(25, 0, 6, 1, 8, 7)_{\adfsmFJh}$,
$(15, 14, 23, 21, 18, 13)_{\adfsmFJh}$,
$(3, 23, 4, 11, 7, 10)_{\adfsmFJh}$,\adfsplit
$(16, 13, 7, 5, 20, 22)_{\adfsmFJh}$,
$(2, 26, 0, 11, 14, 13)_{\adfsmFJh}$,
$(1, 21, 4, 5, 8, 18)_{\adfsmFJh}$,


\adfLgap 
$(21, 11, 3, 5, 20, 12)_{\adfsmFJj}$,
$(13, 9, 24, 12, 0, 1)_{\adfsmFJj}$,
$(25, 7, 14, 0, 16, 1)_{\adfsmFJj}$,\adfsplit
$(11, 12, 14, 15, 20, 17)_{\adfsmFJj}$,
$(4, 10, 17, 26, 19, 21)_{\adfsmFJj}$,
$(2, 1, 18, 23, 10, 20)_{\adfsmFJj}$,

\adfLgap \noindent by the mapping:
$x \mapsto x + 4 j \adfmod{20}$ for $x < 20$,
$x \mapsto (x +  j \adfmod{5}) + 20$ for $20 \le x < 25$,
$x \mapsto x$ for $x \ge 25$,
$0 \le j < 5$.

\ADFvfyParStart{(27, ((6, 5, ((20, 4), (5, 1), (2, 2)))), ((4, 5), (7, 1)))} 

\adfDgap
\noindent{\boldmath $ K_{4^{6},5}$}~
With the point set $Z_{29}$
partitioned into
 residue classes modulo $6$ for $\{0, 1, \dots, 23\}$, and
 $\{24, 25, 26, 27, 28\}$,
 the designs are generated from

\adfLgap 
$(5, 10, 25, 1, 8, 3)_{\adfsmFJc}$,
$(19, 27, 18, 23, 8, 11)_{\adfsmFJc}$,
$(13, 26, 0, 23, 4, 5)_{\adfsmFJc}$,\adfsplit
$(18, 26, 1, 3, 14, 15)_{\adfsmFJc}$,
$(10, 12, 20, 3, 13, 2)_{\adfsmFJc}$,
$(3, 13, 28, 8, 14, 0)_{\adfsmFJc}$,


\adfLgap 
$(12, 26, 6, 9, 21, 22)_{\adfsmFJf}$,
$(5, 27, 6, 4, 3, 15)_{\adfsmFJf}$,
$(2, 25, 5, 6, 12, 13)_{\adfsmFJf}$,\adfsplit
$(3, 19, 16, 0, 12, 14)_{\adfsmFJf}$,
$(11, 26, 5, 0, 7, 10)_{\adfsmFJf}$,
$(3, 28, 17, 4, 10, 13)_{\adfsmFJf}$,


\adfLgap 
$(24, 12, 11, 2, 7, 23)_{\adfsmFJj}$,
$(26, 9, 4, 2, 1, 5)_{\adfsmFJj}$,
$(20, 9, 6, 23, 22, 5)_{\adfsmFJj}$,\adfsplit
$(8, 15, 27, 4, 5, 16)_{\adfsmFJj}$,
$(1, 8, 23, 28, 10, 5)_{\adfsmFJj}$,
$(3, 1, 6, 27, 10, 11)_{\adfsmFJj}$,

\adfLgap \noindent by the mapping:
$x \mapsto x + 4 j \adfmod{24}$ for $x < 24$,
$x \mapsto (x + 2 j \adfmod{4}) + 24$ for $24 \le x < 28$,
$28 \mapsto 28$,
$0 \le j < 6$.

\ADFvfyParStart{(29, ((6, 6, ((24, 4), (4, 2), (1, 1)))), ((4, 6), (5, 1)))} 

\adfDgap
\noindent{\boldmath $ K_{4^{6},5}$}~
With the point set $Z_{29}$
partitioned into
 residue classes modulo $6$ for $\{0, 1, \dots, 23\}$, and
 $\{24, 25, 26, 27, 28\}$,
 the design is generated from

\adfLgap 
$(27, 23, 18, 13, 20, 15)_{\adfsmFJh}$,
$(20, 23, 4, 10, 26, 19)_{\adfsmFJh}$,
$(1, 11, 10, 14, 25, 18)_{\adfsmFJh}$,

\adfLgap \noindent by the mapping:
$x \mapsto x + 2 j \adfmod{24}$ for $x < 24$,
$x \mapsto (x +  j \adfmod{3}) + 24$ for $24 \le x < 27$,
$x \mapsto (x - 27 +  j \adfmod{2}) + 27$ for $x \ge 27$,
$0 \le j < 12$.

\ADFvfyParStart{(29, ((3, 12, ((24, 2), (3, 1), (2, 1)))), ((4, 6), (5, 1)))} 

\adfDgap
\noindent{\boldmath $ K_{4^{6},10}$}~
With the point set $Z_{34}$
partitioned into
 residue classes modulo $6$ for $\{0, 1, \dots, 23\}$, and
 $\{24, 25, \dots, 33\}$,
 the design is generated from

\adfLgap 
$(4, 11, 8, 3, 6, 7)_{\adfsmFJm}$,
$(23, 24, 8, 0, 10, 26)_{\adfsmFJm}$,
$(23, 3, 13, 12, 30, 28)_{\adfsmFJm}$,\adfsplit
$(17, 14, 29, 32, 0, 22)_{\adfsmFJm}$,
$(16, 21, 0, 26, 27, 7)_{\adfsmFJm}$,
$(1, 12, 17, 15, 33, 25)_{\adfsmFJm}$,

\adfLgap \noindent by the mapping:
$x \mapsto x + 3 j \adfmod{24}$ for $x < 24$,
$x \mapsto (x - 24 + 5 j \adfmod{10}) + 24$ for $x \ge 24$,
$0 \le j < 8$.

\ADFvfyParStart{(34, ((6, 8, ((24, 3), (10, 5)))), ((4, 6), (10, 1)))} 

\adfDgap
\noindent{\boldmath $ K_{4^{6},15}$}~
With the point set $Z_{39}$
partitioned into
 residue classes modulo $6$ for $\{0, 1, \dots, 23\}$, and
 $\{24, 25, \dots, 38\}$,
 the design is generated from

\adfLgap 
$(0, 8, 36, 31, 21, 11)_{\adfsmFJm}$,
$(7, 34, 9, 22, 17, 37)_{\adfsmFJm}$,
$(2, 16, 17, 30, 24, 21)_{\adfsmFJm}$,\adfsplit
$(12, 17, 16, 30, 28, 19)_{\adfsmFJm}$,
$(0, 1, 17, 2, 38, 37)_{\adfsmFJm}$,

\adfLgap \noindent by the mapping:
$x \mapsto x + 2 j \adfmod{24}$ for $x < 24$,
$x \mapsto (x - 24 + 5 j \adfmod{15}) + 24$ for $x \ge 24$,
$0 \le j < 12$.

\ADFvfyParStart{(39, ((5, 12, ((24, 2), (15, 5)))), ((4, 6), (15, 1)))} 

\adfDgap
\noindent{\boldmath $ K_{1^{39},21}$}~
With the point set $Z_{60}$
partitioned into
 residue classes modulo $39$ for $\{0, 1, \dots, 38\}$, and
 $\{39, 40, \dots, 59\}$,
 the design is generated from

\adfLgap 
$(15, 20, 46, 57, 2, 1)_{\adfsmFJm}$,
$(11, 52, 49, 18, 34, 3)_{\adfsmFJm}$,
$(7, 48, 10, 5, 6, 16)_{\adfsmFJm}$,\adfsplit
$(0, 10, 20, 22, 40, 44)_{\adfsmFJm}$,

\adfLgap \noindent by the mapping:
$x \mapsto x +  j \adfmod{39}$ for $x < 39$,
$x \mapsto (x - 39 + 7 j \adfmod{21}) + 39$ for $x \ge 39$,
$0 \le j < 39$.

\ADFvfyParStart{(60, ((4, 39, ((39, 1), (21, 7)))), ((1, 39), (21, 1)))} 

\adfDgap
\noindent{\boldmath $ K_{1^{55},25}$}~
With the point set $Z_{80}$
partitioned into
 residue classes modulo $55$ for $\{0, 1, \dots, 54\}$, and
 $\{55, 56, \dots, 79\}$,
 the design is generated from

\adfLgap 
$(37, 67, 39, 2, 47, 59)_{\adfsmFJm}$,
$(38, 68, 40, 3, 48, 60)_{\adfsmFJm}$,
$(39, 69, 41, 4, 49, 61)_{\adfsmFJm}$,\adfsplit
$(40, 70, 42, 5, 50, 62)_{\adfsmFJm}$,
$(41, 71, 43, 6, 51, 63)_{\adfsmFJm}$,
$(74, 41, 36, 16, 40, 13)_{\adfsmFJm}$,\adfsplit
$(75, 42, 37, 17, 41, 14)_{\adfsmFJm}$,
$(76, 43, 38, 18, 42, 15)_{\adfsmFJm}$,
$(66, 44, 39, 19, 43, 16)_{\adfsmFJm}$,\adfsplit
$(67, 45, 40, 20, 44, 17)_{\adfsmFJm}$,
$(37, 64, 51, 6, 8, 46)_{\adfsmFJm}$,
$(38, 65, 52, 7, 9, 47)_{\adfsmFJm}$,\adfsplit
$(39, 55, 53, 8, 10, 48)_{\adfsmFJm}$,
$(40, 56, 54, 9, 11, 49)_{\adfsmFJm}$,
$(41, 57, 0, 10, 12, 50)_{\adfsmFJm}$,\adfsplit
$(59, 42, 41, 0, 5, 7)_{\adfsmFJm}$,
$(60, 43, 42, 1, 6, 8)_{\adfsmFJm}$,
$(61, 44, 43, 2, 7, 9)_{\adfsmFJm}$,\adfsplit
$(62, 45, 44, 3, 8, 10)_{\adfsmFJm}$,
$(63, 46, 45, 4, 9, 11)_{\adfsmFJm}$,
$(78, 34, 46, 12, 40, 18)_{\adfsmFJm}$,\adfsplit
$(77, 46, 35, 52, 19, 13)_{\adfsmFJm}$,
$(79, 52, 53, 25, 36, 14)_{\adfsmFJm}$,
$(68, 39, 11, 33, 22, 50)_{\adfsmFJm}$,\adfsplit
$(66, 7, 29, 40, 18, 51)_{\adfsmFJm}$,
$(66, 9, 20, 26, 37, 48)_{\adfsmFJm}$,

\adfLgap \noindent by the mapping:
$x \mapsto x + 5 j \adfmod{55}$ for $x < 55$,
$x \mapsto (x + 5 j \adfmod{11}) + 55$ for $55 \le x < 66$,
$x \mapsto (x + 5 j \adfmod{11}) + 66$ for $66 \le x < 77$,
$x \mapsto x$ for $x \ge 77$,
$0 \le j < 11$.

\ADFvfyParStart{(80, ((26, 11, ((55, 5), (11, 5), (11, 5), (3, 3)))), ((1, 55), (25, 1)))} 

\adfDgap
\noindent{\boldmath $ K_{1^{99},41}$}~
With the point set $Z_{140}$
partitioned into
 residue classes modulo $99$ for $\{0, 1, \dots, 98\}$, and
 $\{99, 100, \dots, 139\}$,
 the design is generated from

\adfLgap 
$(57, 67, 123, 26, 7, 44)_{\adfsmFJm}$,
$(17, 61, 99, 115, 91, 54)_{\adfsmFJm}$,
$(101, 14, 57, 50, 37, 42)_{\adfsmFJm}$,\adfsplit
$(66, 34, 128, 68, 23, 61)_{\adfsmFJm}$,
$(81, 55, 116, 129, 46, 80)_{\adfsmFJm}$,
$(79, 51, 120, 127, 96, 50)_{\adfsmFJm}$,\adfsplit
$(54, 32, 16, 8, 131, 122)_{\adfsmFJm}$,
$(77, 98, 29, 11, 71, 25)_{\adfsmFJm}$,
$(52, 124, 55, 12, 71, 102)_{\adfsmFJm}$,

\adfLgap \noindent by the mapping:
$x \mapsto x +  j \adfmod{99}$ for $x < 99$,
$x \mapsto (x +  j \adfmod{11}) + 99$ for $99 \le x < 110$,
$x \mapsto (x - 110 + 10 j \adfmod{30}) + 110$ for $x \ge 110$,
$0 \le j < 99$.

\ADFvfyParStart{(140, ((9, 99, ((99, 1), (11, 1), (30, 10)))), ((1, 99), (41, 1)))} 



\end{document}